\documentclass[10pt]{article}

\usepackage{amsmath}
\usepackage{amssymb}
\usepackage{indentfirst}
\usepackage{graphics} 
\usepackage{color}

\makeatletter

\@addtoreset{equation}{section}
\makeatother
\def\R{\mbox{\boldmath $R$}}

\def\<{\langle}
\def\>{\rangle}

\newtheorem{lem}{Lemma}[section]
\newtheorem{theo}{Theorem}[section]
\newtheorem{rem}{Remark}[section]
\newtheorem{pro}{Proposition}[section]

\makeatletter
   
   \@addtoreset{equation}{section}
\makeatother

\begin{document}
\title{\bf  A dissipative logarithmic-Laplacian type of plate equation: asymptotic profile and decay rates}

\author{Ruy Coimbra Char\~ao\thanks{Corresponding author: ruy.charao@ufsc.br} \; and \; Alessandra Piske\thanks {alessandrapiske@gmail.com}  \\{\small Department of Mathematics} \\{\small Federal University of Santa Catarina} \\ {\small 88040-270, Florianopolis, Brazil,} 
\\
and\\Ryo Ikehata\thanks{ikehatar@hiroshima-u.ac.jp} \\ {\small Department of Mathematics}\\ {\small Division of Educational Sciences}\\ {\small Graduate School of Humanities and Social Sciences} \\ {\small Hiroshima University} \\ {\small Higashi-Hiroshima 739-8524, Japan}}
\date{}
\maketitle
\begin{abstract}
We introduce a new model of the logarithmic type of wave like plate equation with a nonlocal logarithmic damping mechanism. We consider the Cauchy problem for this new model in ${\bf R}^{n}$, and study the asymptotic profile and optimal decay rates of solutions as $t \to \infty$ in $L^{2}$-sense. The operator $L$ considered in this paper was first introduced to dissipate the solutions of the wave equation in the paper studied by Char\~ao-Ikehata \cite{Log-damping}. We will discuss the asymptotic property of the solution as time goes to infinity to our Cauchy problem, and in particular, we classify the property of the solutions into three parts from the viewpoint of regularity of the initial data, that is,  diffusion-like, wave-like, and both of them.  

\end{abstract}
\section{Introduction}
\footnote[0]{Keywords and Phrases: Wave-like equation; Logarithmic damping; $L^{2}$-decay; asymptotic profile, optimal estimates, regularity.}
\footnote[0]{2010 Mathematics Subject Classification. Primary 35L05; Secondary 35B40, 35C20, 35S05.}

We consider in this work a new model of evolution equations based on an operator $L$, that combines the composition of logarithm function with the Laplace operator as follows:
\begin{align}
& u_{tt} + Lu + (I+L)^{-1}u_{t} = 0,\ \ \ (t,x)\in (0,\infty)\times {\bf R}^{n},\label{eqn}\\
& u(0,x)= u_{0}(x), \quad  u_{t}(0,x)= u_{1}(x),\ \ \ x\in{\bf R}^{n},\label{initial}
\end{align}
where the linear operator  
\[L: D(L) \subset L^{2}({\bf R}^{n}) \to L^{2}({\bf R}^{n})\]
is defined by 
\[D(L) := \left\{f \in L^{2}({\bf R}^{n}) \,\bigm|\,\int_{{\bf R}^{n}}(\log(1+\vert\xi\vert^{2}))^{2}\vert\hat{f}(\xi)\vert^{2}d\xi < +\infty\right\},\]
and for $f \in D(L)$,  
\[(Lf) (x) := {\cal F}_{\xi\to x}^{-1}\left(\log (1+\vert\xi\vert^{2})\hat{f}(\xi)\right)(x).\]
\noindent
Here, one has just denoted the Fourier transform ${\cal F}_{x\to\xi}(f)(\xi)$ of $f(x)$ by 
\[{\cal F}_{x\to\xi}(f)(\xi) = \hat{f}(\xi) := \displaystyle{\int_{{\bf R}^{n}}}e^{-ix\cdot\xi}f(x)dx, \quad \xi \in {\bf R}^n,\]
as usual with $i := \sqrt{-1}$, and ${\cal F}_{\xi\to x}^{-1}$ expresses its inverse Fourier transform. 
\noindent
Since the operator $L$ is non-negative and self-adjoint in $L^{2}({\bf R}^{n})$ (see \cite{Log-damping}), the square root 
$$L^{1/2}: D(L^{1/2}) \subset L^{2}({\bf R}^{n}) \to L^{2}({\bf R}^{n})$$
can be defined, and is also nonnegative and self-adjoint with its domain
\[D(L^{1/2}) = \left\{f \in L^{2}({\bf R}^{n}) \,\bigm|\,\int_{{\bf R}^{n}}\log(1+\vert\xi\vert^{2})\vert\hat{f}(\xi)\vert^{2}d\xi < +\infty\right\}.\] 
Note that $D(L^{1/2})$ becomes Hilbert space with its graph norm
$$\Vert v\Vert_{D(L^{1/2})} := \left(\Vert v\Vert^{2} + \Vert L^{1/2}v\Vert^{2}\right)^{1/2},$$
where the $L^{2}({\bf R}^{n})$-norm is defined by
\[\Vert\cdot\Vert := \Vert\cdot\Vert_{L^{2}({\bf R})}.\]

It is easy to check that 
$$H^{s}({\bf R}^{n}) \hookrightarrow D(L^{1/2}) \hookrightarrow  L^{2}({\bf R}^{n})$$
for $s > 0$.\\
\noindent
Symbolically writing, one can see
\[L = \log(I-\Delta),\]
where $\Delta$ is the usual Laplace operator defined on $H^2({\bf R}^n)$.

As for the existence of the unique solution to problem \eqref{eqn}-\eqref{initial} we discuss, on the next section, by employing the standard Lumer-Phillips Theorem. 

The associated energy inequality to the system  \eqref{eqn}-\eqref{initial}  is
\begin{equation}\label{energy}
E_{u}(t) \leq E_{u}(0),
\end{equation} 
where
\[
E_{u}(t) := \frac{1}{2}\left(\Vert u_{t}(t,\cdot)\Vert_{L^{2}}^{2} + \Vert \log^{1/2}(I-\Delta) u(t,\cdot)\Vert_{L^{2}}^{2}\right).
\]
The inequality \eqref{energy} implies that the the total energy is a non increasing function in time because of the existence of some kind of dissipative term $ (I+L)^{-1}u_{t}$.  

A main topic of this paper is to find an asymptotic profile of the solution in the $L^{2}$ framework to problem \eqref{eqn}-\eqref{initial}, and to apply it to investigate the optimal rate of decay of the solution in terms of the $L^{2}$-norm. We study the equation \eqref{eqn} only from the purely mathematical point of view. The model equation itself is strongly inspired from the related paper due to Dharmawardane-Nakamura-Kawashima \cite{DNK}.\\  

A motivation of this research has its origin in the study of the damped plate equation under effects  of rotational inertia
\begin{equation*}\label{strong}
u_{tt}-\Delta u_{tt}-\Delta u + \Delta^2 u +  u_{t} = 0 
\end{equation*}
that can be written as 
\begin{equation}\label{strong}
u_{tt}-\Delta u + (I-\Delta)^{-1} u_{t} = 0.
\end{equation}
In a work Char\~ao-Espinoza-Ikehata \cite{C-E-I} they study a more general model than \eqref{strong} with a super damping. A pioneering work on the dissipative structure and nonlinear property of the more generalized systems including \eqref{strong} is discussed in Dharmawardane \cite{D}, which is mentioning a regularity loss-structure of the equation. After \cite{D}, Fukushima-Ikehata-Michihisa \cite{FIM-1} have studied the asymptotic profiles of the solution to \eqref{strong}, and such profiles are divided into two parts: one is the Gauss kernel like for high regularity initial data, and the other is related with the non-diffusive oscillating property for low regularity initial data from the viewpoint of regularity-loss structure. In this connection, such a regularity-loss structure has been first discovered and named by S. Kawashima through the analysis for the dissipative Timoshenko system.

An analysis of the dissipative mechanism of \eqref{strong} goes back to the two pioneering works of G. Ponce \cite{Po} and Y. Shibata \cite{S}, where they studied various $L^{p}$-$L^{q}$ estimates of the solution to the Cauchy problem for the equation:
\begin{equation}\label{strong-0}
u_{tt}-\Delta u -\Delta u_{t} = 0.
\end{equation}
 After them, an asymptotic profile and the optimal estimates of the solution can be introduced in the papers \cite{ITY}, \cite{I-14} and \cite{IO}. They investigated a singularity near $0$-frequency region of the solution to \eqref{strong-0} in terms of $L^{2}$-norm of solutions. In this connection, in \cite{ BV-1, BV-2} and \cite{Mi} a higher order asymptotic expansion of the solution as $t \to \infty$ to the equation \eqref{strong-0} is precisely investigated.

On the other hand, the so-called critical exponent problem for semi-linear equations of \eqref{strong-0} is first developed by D'Abbicco-Reissig \cite{DR}, and this paper has been the beginning of a series of related papers studying structurally damped wave models with nonlinearity. Unfortunately, at present nobody knows the precise value of the critical exponent $p^{*}$ of the equation \eqref{strong-0} with power type nonlinearity $\vert u\vert^{p}$. This study is based on the $L^{p}-L^{q}$-estimates derived in \cite{S}.

Recently, the equation \eqref{strong-0} is generalized to the linear and semi-linear models, respectively:
\begin{equation}\label{strong-2}
u_{tt} + (-\Delta)^{\sigma}u + (-\Delta)^{\theta}u_{t} = 0,
\end{equation}
\begin{equation}\label{strong-3}
u_{tt} + (-\Delta)^{\sigma}u + (-\Delta)^{\theta}u_{t} = f(u,u_{t}).
\end{equation}
A study on asymptotic profile and $L^{p}$-$L^{q}$ estimates to the equation \eqref{strong-2} has been done in the papers \cite{CLI}, \cite{DE}, \cite{DEP}, \cite{DGL}, \cite{LIC}, \cite{NR}, and \cite{IT}, and the corresponding critical exponent problems (mainly) to the equation \eqref{strong-3} are treated in the papers \cite{DE-2}, \cite{AR, AR-2}, \cite{K}, and \cite{PMR}. As for the related motivated topics concerning the asymptotic profiles for the higher order PDEs and the other types of linear equations one can also cite \cite{C-1,C-2}, \cite{R}, \cite{TY}, \cite{SW}, and the references therein. 

In \cite{II} and \cite{FIM}, the so-called regularity-loss structure of the solution in the high frequency zone can be studied to \eqref{strong-2} with $\sigma = 1$ and $\theta > 1$, and these researches are strongly inspired from the abstract theory due to \cite{GGH}.

Quite recently Char\~ao-Ikehata  \cite{Log-damping} introduced a new type of damping term of logarithm type to the wave equation, and it is expressed in the Fourier space as follows:
\begin{equation}\label{logdamping}
\hat{u}_{tt} + \vert\xi\vert^{2}\hat{u} + \log(1+\vert\xi\vert^{2})\hat{u}_{t} = 0.
\end{equation}
Symbolically writing, one sees
\[u_{tt}-\Delta u + \log(I-\Delta)u_{t} = 0.\]
In \cite{CDI}, \eqref{logdamping} is more generalized to the equation such that
\[u_{tt}-\Delta u + \log(I + (-\Delta)^{\theta})u_{t} = 0\]
for $\theta > 1/2$.  

On reconsidering our problem \eqref{eqn}-\eqref{initial} in the Fourier space, our equation becomes  
\begin{equation}\label{log-logdamping}
\hat{u}_{tt} + \log(1+\vert\xi\vert^{2})\hat{u} + (1+\log(1+\vert\xi\vert^{2}))^{-1}\hat{u}_{t} = 0.
\end{equation} 
We should  investigate characteristics roots of \eqref{log-logdamping}  to see where they are complex-valued or not,   for small $\xi \in {\bf R}^{n}$  and on  large frequency zone.


In order to introduce our main results we should define function spaces with respect to the logarithmic Laplace operator$L$ such that for $\delta \geq 0$
$$Y^{\delta}  = \{ f \in L^2({\bf R}^n); \int _{{\bf R}_{\xi}^{n}} (1+\log (1+|\xi|^2 ))^{\delta} |\hat{f} (\xi)|^2d \xi < \infty  \} $$
with its natural norm
\[\Vert f\Vert_{Y^{\delta}} := \left(\int _{{\bf R}_{\xi}^{n}} (1+\log (1+|\xi|^2 ))^{\delta} |\hat{f} (\xi)|^2d \xi\right)^{1/2}\]
for $f \in Y^{\delta}$. 
\begin{rem}{\rm Due to the fact that  $ \log (1+|\xi |^2) \leq |\xi|^2  $ for all $\xi \in {\bf R}^n$,  one notices $ H^\delta ({\bf R}^{n}) \subset Y^\delta \subset L^2$ for $\delta \geq 0$.}
\end{rem}
\noindent
Furthermore, we set
\[I_{0,l} :=  \sqrt{\| u_0 \|_{1,1}^2 + \| u_1 \|_{1,1}^2 + \|u_0 \|_{Y^{l+1}}^2 + \| u_1 \|_{Y^l}^2 }\quad (l \geq 0),\]
and
\[P_{j} := \int_{{\bf R}^{n}}u_{j}(x)dx\quad (j = 0,1).\]

Our results read as follows.
\begin{theo}\label{1.1} Let $ n \geq 1 $ and $l \geq 1$. If $ (u_0, u_1) \in (L^{1,1}({\bf R}^{n}) \cap Y^{l+1}) \times (L^{1,1}({\bf R}^{n}) \cap Y^{l})  $, then there exists a constant $C>0 $, which is independent of $t, u_0, u_1$ such that 
\[\|u(t,\cdot) - {\cal F}_{\xi\rightarrow x}^{-1}(\varphi _1(t,\xi))(\cdot)\|_2 \]
\[\leq \left\{\begin{matrix}
CI_{0,l}^2 t^{-\frac{n+2}{4}} & \text{ if } l\geq 1 \text{ and } n\leq 2; \text{ if } n \geq 3 \text{ and } l \geq n/2,     \\ 
CI_{0,l}^2 t^{-\frac{l+1}{2}} &  \text{ if } n \geq 4 \text{ and } n/2 -1 < l \leq n/2; \text{  if } n=3 \text{ and } 1 \leq l \leq 3/2,
\end{matrix}\right.\]
for $t \gg 1$, where
\[\varphi _1 (t,\xi) := e^{-t\log(1+|\xi|^2)(1+\log(1+|\xi |^2))} (P_0 + P_1).\] 
\end{theo} 

\begin{theo}\label{1.2} Let $ n \geq 5 $ and $l \geq 1$. If $ (u_0, u_1) \in (L^{1,1}({\bf R}^{n}) \cap Y^{l+1}) \times (L^{1,1}({\bf R}^{n}) \cap Y^{l})  $, then there exists a constant $C>0 $, which is independent of $t, u_0, u_1$ such that 
\[\| u(t,\cdot) - {\cal F}_{\xi\rightarrow x}^{-1}(\varphi _2(t,\xi))(\cdot) \|_2\]
\[\leq \left\{\begin{matrix}
CI_{0,l}^2 t^{-\frac{n}{4}} & \text{if }   1\leq l < n/2-1 \text{ and } 5\leq n \leq 8; \text{ if }  n>8  \text{ and } n/2-3 < l < n/2-1, \\ 
CI_{0,l}^2 t^{-\frac{l+3}{2}} &  \text{ if }   n>8  \text{ and } 1\leq l \leq n/2-3,
 \end{matrix}\right.\]
for $t \gg 1$, where
\[\varphi _2(t, \xi) :=e^{-\frac{t}{2 \log (1+|\xi |^2 )}} \frac{\sin (\sqrt{\log (1+|\xi |^2 )}t )}{\sqrt{\log (1+|\xi |^2)}} \hat{u}_1(\xi) +  e^{-\frac{t}{2 \log (1+|\xi |^2 )}} \cos (\sqrt{\log (1+|\xi |^2 )}t ) \hat{u}_0(\xi).\] 
\end{theo} 

\begin{theo}\label{1.3} Let $ n \geq 4 $ and $ l= \frac{n}{2}-1 $. If $ (u_0, u_1) \in (L^{1,1}({\bf R}^{n}) \cap Y^{l+1}) \times (L^{1,1}({\bf R}^{n}) \cap Y^{l})  $, then there exists a constant $C>0 $, which is independent of $t, u_0, u_1$ such that 
\begin{align*}
\| u(t, \cdot ) - {\cal F}_{\xi\rightarrow x}^{-1}(\varphi (t,\xi))(\cdot) \|_2 \leq CI_{0,l}^2 t^{-\frac{n+2}{4}}
\end{align*}
for $t \gg 1$, where
\[\varphi (t,\xi) := \varphi_{1}(t,\xi) + \varphi_{2}(t,\xi).\] 
\end{theo} 
\begin{rem} {\rm In Theorems above one has assumed $l \geq 1$ because in this case the problem \eqref{eqn}-\eqref{initial} has a unique weak solution in the class 
\[u \in C([0,+\infty),Y^{2}) \cap C^{1}([0,+\infty),Y^{1}) \cap C^{2}([0,+\infty),L^{2}({\bf R}^{n})).\]
For details concerning the unique existence of solutions, we will discuss in Theorem 2.1 in subsection 2.1.}
\end{rem}
\begin{rem} {\rm The value $l^{*}(n)$ defined by $l^{*}(n) := \frac{n}{2}-1$ expresses a kind of critical number on the regularity $l \geq 1$, which divides the property of the solution $u(t,x)$ into three types: one is diffusive-like (Theorem \ref{1.1}), the other is wave-like (Theorem \ref{1.2}), and the remaining is both of them (Theorem \ref{1.3}). }
\end{rem}
Throughout the discussions to get Theorems \ref{1.1}, \ref{1.2} and \ref{1.3}, one can obtain the following crucial results, which shows the optimality concerning the $L^{2}$-decay rates of the solution $u(t,x)$ to problem (1.1)-(1.2).  
\begin{theo}\label{TEO0} Let $ n \geq 4 $ and $ 1 \leq l \leq \displaystyle{\frac{n}{2}}-1 $. If $ (u_0, u_1) \in (L^{1,1}({\bf R}^{n}) \cap Y^{l+1}) \times (L^{1,1}({\bf R}^{n}) \cap Y^{l})  $, then the solution $u(t,x)$ to problem {\rm (1.1)-(1.2)}  satisfies 
\begin{align*}
\| u(t,\cdot )\|_2 \leq C I_{0,l} t^{-\frac{l+1}{2}}
\end{align*}
for $t \gg 1 $, where $ C $ is a positive constant which depends only on $n$.  
\end{theo}
\begin{rem}{\rm The rate of decay obtained in Theorem \ref{TEO0} seems exactly optimal, however, one cannot obtain the lower bound of time-decay rate. This is still open.}
\end{rem}
\begin{theo}\label{TEO1} Let $ 1\leq n \leq 3 $. If $ (u_0, u_1) \in (L^{1,1}({\bf R}^{n}) \cap Y^{2}) \times (L^{1,1}({\bf R}^{n}) \cap Y^{1})  $, then there exists  constants $C_1,C_2>0 $ independent of $t$ such that 
\begin{align*}
C_1 |P_0+P_1| t^{-\frac{n}{4}}  \leq \| u(t, \cdot ) \|_2 \leq C_2 I_{0,1}  t^{-\frac{n}{4}} 
\end{align*}
for all $t \gg 1$ provided that $P_{1}+P_{0} \ne 0$.
\end{theo}

\begin{theo}\label{TEO2} Let $ n \geq 4 $ and $ \varepsilon >0 $. If $ (u_0, u_1) \in (L^{1,1}({\bf R}^{n}) \cap Y^{\frac{n}{2}+\varepsilon}) \times (L^{1,1}({\bf R}^{n}) \cap Y^{\frac{n-2}{2}+\varepsilon})  $, then there exists  constants $C_1,C_2>0 $ independent of $t$ such that 
\begin{align*}
C_1 |P_0+P_1| t^{-\frac{n}{4}}  \leq \| u(t, \cdot ) \|_2 \leq C_2 I_{0,\frac{n}{2} +\varepsilon-1}  t^{-\frac{n}{4}} 
\end{align*}
for all $t \gg 1$,  provided that $P_{1}+P_{0} \ne 0$.
\end{theo}


This paper is organized as follows. In Section 2 we prepare several important lemmas, which will be used later, and in particular, these lemmas are closely related with hypergeometric functions (see \cite{Log-damping}). In the final part of Section 2 we discuss the unique existence of solutions to problem \eqref{eqn}-\eqref{initial} using ideas from Luyo Sanchez \cite{Luyo}. 
In Section 3 we obtain decay estimates for the  $L^2$-norm of solutions by using the multiplier method in the Fourier space (cf. \cite{UKS}) only for $n\geq 3$ in the case when the initial data have a low regularity. The asymptotic profile and related estimates of solution are obtained in Section 4, and in particular, Theorems \ref{1.1}, \ref{1.2} and \ref{1.3} are proved in subsection 4.4, and in the final subsection 4.5 we give the proof of Theorems \ref{TEO0}, \ref{TEO1} and \ref{TEO2}. 
\vspace{0.2cm}
\par
{\bf Notation.} {\small Throughout this paper, $\| \cdot\|_q$ stands for the usual $L^q({\bf R}^{n})$-norm. For simplicity of notation, in particular, we use $\| \cdot\|$ instead of $\| \cdot\|_2$. We also introduce the following weighted functional spaces:
\[L^{1,\gamma}({\bf R}^{n}) := \left\{f \in L^{1}({\bf R}^{n}) \; \bigm| \; \Vert f\Vert_{1,\gamma} := \int_{{\bf R}^{n}}(1+\vert x\vert^{\gamma})\vert f(x)\vert dx < +\infty\right\}.\]
Furthermore, we denote the surface area of the $n$-dimensional unit ball by $\omega_{n} := \displaystyle{\int_{\vert\omega\vert = 1}}d\omega$. 

}


\section{Basic results and existence of solutions}
In this section we shall collect important lemmas to derive precise estimates of various quantities related to the solution to problem \eqref{eqn}-\eqref{initial}. These are already studied and developed in our previous works (see \cite{Log-damping, CDI}).

The following main estimate to the function  
$$I_p(t)= \int_0^{1}(1+r^{2})^{-t}r^p dr$$ 
is a direct consequence of the cases $p \geq 0$ in Char\~ao-Ikehata \cite{Log-damping} and $-1<p<0$ in Char\~ao-D'Abbicco-Ikehata \cite{CDI}. In the following  notation $A \approx B$ means that  $c_1 A \leq B \leq c_2 A$ for some positive constants $c_1, c_2$.

\begin{lem}\label{general-p}
 Let $p > -1$ be a real number. Then it holds that
$$I_p(t) \approx t^{-\frac{p+1}{2}}, \quad t \gg 1.$$
\end{lem}

In order to deal with the high frequency part of estimates, one relies on the function again 
$$J_p(t)=\int_1^{\infty}(1+r^2)^{-t}r^p dr$$
for $p \in {\bf R}$.

Then the next lemma is important to get estimates on the zone of high frequency to problem \eqref{eqn}--\eqref{initial}. The proof  appears in Char\~ao-Ikehata \cite{Log-damping}.
\begin{lem}\label{infit}
\,Let $p \in {\bf R}$. Then it holds that 
$$J_p(t) \approx \dfrac{2^{-t}}{t-1}, \quad t \gg 1.$$
\end{lem}
\vspace{0.2cm}
For later use we prepare the following simple lemma, which implies the exponential decay of the middle frequency part.
\begin{lem}\label{intermid}\,Let $p \in {\bf R}$, and $\eta \in (0,1]$. Then there is a constant $C > 0$ such that 
$$\int_{\eta}^{1}(1+r^{2})^{-t}r^{p}dr \leq C(1+\eta^{2})^{-t}, \quad t \geq 0.$$
\end{lem}

\begin{rem}
{\rm We note that the proof of Lemma \ref{general-p} is done by using simple differential calculus and the theory from hypergeometric functions (see Watson \cite{W}).}

\end{rem}

\begin{lem}\label{LemmaTecnico1}
Let $ c , \nu  $ be positive real numbers and $a \in {\bf R}$. Then, there exists a constant  $ C>0 $ such that 
$$ t ^{\nu} e ^{-c (1+\log (1+ |\xi |^2  ))^{a}t} \leq C (1+ \log (1+ | \xi | ^2 ))^{-a \nu}.$$ 
\end{lem}
{\it Proof.}\
We set $ s:= c (1+\log (1+ |\xi |^2  ))^{a}t $. Then $ t = c^{-1} (1+ \log (1+ |\xi |^2  ))^{-a} s$ and 
$$ t ^{\nu } = c^{-\nu} (1+ \log (1+ |\xi |^2  ))^{-a \nu } s^{\nu} . $$

The definition of $s$ implies 
$$  t ^{\nu } e ^{-c (1+\log (1+ |\xi |^2  ))^{a}t} = c^{-\nu} (1+ \log (1+ |\xi |^2  ))^{-a \nu } s^{\nu} e^{-s} .$$
Since the function  $ {\bf R} \ni s \mapsto s^{\nu} e^{-s}  $ is bounded, there exists $ C > 0 $ such that 
$$ t ^{\nu } e ^{-c (1+\log (1+ |\xi |^2  ))^{a}t} \leq C (1+ \log (1+ |\xi |^2  ))^{-a \nu } .$$ 
\hfill
$\Box$

\begin{lem}\label{lemmahiperbolicsine}
There exists a constant $C>0 $ such that 
$$ \frac{\sinh x }{x} \leq C e^{x}  $$
for $ x >0 $. 
\end{lem}

Let $f \in L^1({\bf R}^{n})$. Then,  we may decompose the Fourier transform of  $f$ as follows:   
\begin{equation} \label{decompo}
\hat{f}(\xi)=A_f(\xi)-iB_f(\xi)+P_f,
\end{equation}
for all $\xi \in {\bf R}^{n}$, where $i := \sqrt{-1}$ and 
\begin{itemize}
\item[$\bullet$] $A_f(\xi)=\displaystyle{\int_{{\bf R}^{n}}}{(\cos(x\cdot\xi)-1)f(x)}dx,$
\item[$\bullet$] $B_f(\xi)=\displaystyle{\int_{{\bf R}^{n}}}{\sin(x\cdot\xi)f(x)}dx,$
\item[$\bullet$] $P_f=\displaystyle{\int_{{\bf R}^{n}}}{f(x)}dx.$\end{itemize}

We define the weighted $L^{1}$-space $L^{1, \kappa} ({\bf R}^{n}) $ by 
$$ L^{1, \kappa} ({\bf R}^{n}) := \left\{ f \in L^1({\bf R}^{n}) : \int _{{\bf R}^{n}} (1+|x|^{\kappa} ) |f(x)| d x  < +\infty \right\}. $$

The next lemma can be proved in a standard way (see \cite{I-04}).
\begin{lem}\label{lema2.6}
\begin{itemize}
\item[{\rm i)}] If\;  $f \in L^1({\bf R}^{n})$, then for all $\xi \in {\bf R}^{n}$ it is true that  
$$|A_f(\xi)|\leq L\|f\|_{L^1} \quad \text{ and  } \quad |B_f(\xi)|\leq N\|f\|_{L^1}.$$
\item[{\rm ii)}] If \;$0<\kappa \leq 1$ and  $f \in  L^{1,\kappa}({\bf R}^{n})$, then for all  $\xi \in {\bf R}^{n}$ it is true that 
$$|A_f(\xi)|\leq K|\xi|^\kappa\|f\|_{L^{1,\kappa}} \quad \text{ and  } \quad |B_f(\xi)|\leq M|\xi|^\kappa\|f\|_{L^{1,\kappa}}$$
\end{itemize}
\noindent with  $L$, $N$, $K$ and  $M$ positive  constants depending only on the dimension $n$ or $n$ and $\kappa$.\\
\end{lem}


\subsection{Existence and Uniqueness}

We note that the Cauchy problem \eqref{eqn}-\eqref{initial} is equivalent to
\begin{equation}\label{eq2}
\left\{\begin{matrix}
(I+L)u_{tt}+L(I+L)u +u_t = 0, \\ 
u(0,x)= u_0(x),\\ 
u_t(0,x)=u_1(x).
\end{matrix}\right.
\end{equation}
By taking the inner product of the equation in \eqref{eq2} by $u_t$, we obtain
$$ \frac{1}{2} \frac{\mathrm{d} }{\mathrm{d} t} \left ( \| u_t(t,\cdot) \| ^2 + \| L^{1/2}u_t(t,\cdot) \| ^2 + \| Lu(t,\cdot) \|^2 + \|L^{1/2}u(t,\cdot)\|^2   \right ) + \| u_t (t,\cdot)\| ^2 =0 .  $$
We define the total energy as
$$ E(t):= \| u_t(t,\cdot) \| ^2 + \| L^{1/2}u_t(t,\cdot) \| ^2 + \| Lu(t,\cdot) \|^2 + \|L^{1/2}u(t,\cdot)\|^2. $$
Then, we can observe that $E(t)$ is a non-increasing function. 

Associated to \eqref{eq2} one can choose the following energy space
$$ X = Y ^2 \times Y^1 .$$
\noindent
Now, at least formally, from \eqref{eq2} one can write
$$ u _{tt} = -(I+L)^{-1} (L^2+L) u  -(I+L)^{-1}u_t = -(I+L)^{-1} (L^2+L +I) u  -(I+L)^{-1}(u_t-u).    $$
Then if we define $v= u_t $ and $U(t) = \begin{pmatrix}
u\\ v
\end{pmatrix}$ we can reduce the second order equation of \eqref{eq2} to a system of the first order as follows:
\begin{eqnarray*}
\frac{\mathrm{d} U}{\mathrm{d} t} &=&  \begin{pmatrix}
u_t\\ 
-(I+L)^{-1} (L^2+L +I) u  -(I+L)^{-1}(u_t-u) 
\end{pmatrix} \\
&=& \begin{pmatrix}
0 & I\\ 
-(I+L)^{-1}(L^2+L+I) & 0 
\end{pmatrix} \begin{pmatrix}
u\\ 
v
\end{pmatrix} + \begin{pmatrix}
0\\ 
(I+L)^{-1}(u-v)
\end{pmatrix} \\
&=& \begin{pmatrix}
0 & I\\ 
-A & 0 
\end{pmatrix} \begin{pmatrix}
u\\ 
v
\end{pmatrix} + \begin{pmatrix}
0\\ 
(I+L)^{-1}(u-v)
\end{pmatrix} .
\end{eqnarray*}
\noindent
Thus, the first order evolution equation to $U$ can be written as 
\begin{equation} \label{eq-ABJ}
\frac{\mathrm{d} U}{\mathrm{d} t} = B U +JU, \quad U(0)=(u_0, u_1), 
\end{equation}
where  formally the operator $A$ is given by
$$A := (I+L)^{-1}(L^2 + L + I)= L + (I+L)^{-1},$$
and  the operators $ B, J$ are given by
$$ B= \begin{pmatrix}
0 & I\\ -A & 0 
\end{pmatrix} ,  \quad JU =\begin{pmatrix}
0\\ 
(I+L)^{-1}(u-v)
\end{pmatrix} , \quad U \in D(A).$$

We need to give a precise definition of the domain of operator $A$.
To do that we set   
\begin{eqnarray*}
D(A) &=& \{ u\in Y^2: \text{ there exists  } y = y_u \in Y^1 \text{ such that  } \\
& &(Lu, L\psi ) +(L^{1/2}u , L^{1/2} \psi ) +(u, \psi) = (y, \psi) + (L^{1/2}y, L^{1/2}\psi) \text{ } \forall \text{ } \psi \in Y^2 \}.
\end{eqnarray*}
We observe that  $ 0 \in D(A) $, so $ D(A) \neq \phi $ and we may define $$ A : D(A) \rightarrow Y^1, \quad  \mbox{by} \quad  Au = y_u , \; u \in D(A).$$ The following result guarantees us that  $A$ is well defined.  
\begin{lem} For $ u \in Y^2 $, there exists at most  one  $ y \in Y^1 $ that satisfies 
\begin{equation}
(Lu, L\psi ) +(L^{1/2}u , L^{1/2} \psi ) +(u, \psi) = (y, \psi) + (L^{1/2}y, L^{1/2}\psi) \text{ } \forall \text{ } \psi \in Y^2 .
\end{equation}
\end{lem}
{\it Proof.}\
Suppose that $ y_1, y_2 \in Y^1   $ satisfy the above relation. Then $ y= y_1-y_2 $ satisfies
$$ (y, \psi ) + (L^{1/2}y, L^{1/2}\psi )= 0 \text{ \; for each\;\;} \psi \in Y^2.  $$
In particular,  
\begin{equation}\label{eq4}
(y, \psi ) + (L^{1/2}y, L^{1/2}\psi )= 0 \; \text{ for each\;\;} \psi \in C_0^{\infty}({\bf R}^{n}) .
\end{equation}
By density of $ C_0^{\infty}({\bf R}^{n})  $ in $ Y^1 $, there exists a sequence $\{ y_n \} \subset C_0^{\infty}({\bf R}^{n}) $ such that $$ \| y_n - y \|_{Y^1}  \rightarrow 0  .  $$
This implies  $\|y_n\|_{Y^1} \rightarrow \|y\|_{Y^1} $ and $ \| y_n - y \|_{Y^1}^2  \rightarrow 0 $. Thus we  have 
$$ \| y_n \| _{Y^1}^2 - 2 (y_n, y)_{Y^1} + \| y \|_{Y^1} ^2  = \| y_n - y \|_{Y^1}^2 \rightarrow 0,$$
which implies 
$$ \; (y_n, y)_{Y^1} = \| y \|_{Y^1} ^2 . $$

On the other hand, by \eqref{eq4} and the density argument we get
$$ (y,y_n)_{Y^1}=(y, y_n ) + (L^{1/2}y, L^{1/2}y_n )= 0,$$
which implies $\displaystyle{\lim _{n \rightarrow \infty}}(y,y_n)_{Y^1} = 0  $. Therefore, we can conclude  that $ \|y \| _{Y^1} =0 .$
\hfill
$\Box$

\begin{lem}
$ D(A) \subset Y^3 $ and there exists  $ c > 0 $ such that 
$$\| Au \| _{Y^1} \leq c \| u \|_{Y^3} . $$
\end{lem}
{\it Proof.}\
Let $ u \in D(A).  $ Then there is $ y \in Y^1 $ that satisfies 
$$(Lu, L\psi ) +(L^{1/2}u , L^{1/2} \psi ) +(u, \psi) = (y, \psi) + (L^{1/2}y, L^{1/2}\psi) \text{ } \forall \text{ } \psi \in Y^2 .$$

We define $F: Y^1 \rightarrow {\bf R}$ by 
$$ \left \langle F, \psi \right \rangle = (y, \psi ) + (L^{1/2}y, L^{1/2}\psi), \; \psi \in  Y^1. $$ Then $ F $ is well defined, because $y, \psi, L^{1/2}y $ and $L^{1/2}\psi$ are in $L^2({\bf R}^{n}) $ and $F$ is linear. Furthermore, $F$ is a  continuous operator. In  fact   
\begin{eqnarray*}
| \left \langle F, \psi \right \rangle | &\leq & | (y, \psi) | + | (L^{1/2}y, L^{1/2}\psi )| \leq \|y \| \|\psi \| + \| L^{1/2}y\| \| L^{1/2}\psi\|  \\
&=&  \|\hat{y} \| \|\hat{ \psi } \| + \| (\log (1+|\xi|^2))^{1/2}\hat{ y }\| \| (\log (1+|\xi|^2))^{1/2} \hat{ \psi }\| \\
&\leq & 2 \| (1+\log (1+|\xi|^2))^{1/2}\hat{ y }\| \| (1+ \log (1+|\xi|^2))^{1/2} \hat{ \psi }\| \\
&=& 2 \| y \|_{Y^1} \| \psi \| _{Y^1}.
\end{eqnarray*}
\noindent
Since $\mathcal{S}({\bf R}^{n}) \subset  Y^2 \subset Y^1 $, we have $ ( Y^1)' \subset ( Y^2)' \subset \mathcal{S}'({\bf R}^{n})  .$ In other words, $ F $ can be seen as a tempered distribution  and for all $ \psi \in \mathcal{S}({\bf R}^{n})  $ it holds that
 $$ (Lu, L\psi ) +(L^{1/2}u , L^{1/2} \psi ) +(u, \psi)  = \left \langle F, \psi \right \rangle . $$
Thus   
 $$ L^2u+Lu +u  =  F \,\,\, \text{in} \,\,\,\mathcal{S}'({\bf R}^{n} ).  $$
By applying the Fourier transform, the definition of $F$ and the operator $L$, we arrive that 
$$\left [ \log ^2(1+|\xi |^2 )+ \log (1+|\xi |^2 ) + 1  \right ] \hat{u} = \hat{ F} = \left [ 1 + \log (1+|\xi |^2 ) \right ] \hat{ y},$$
that is,  $$  \hat{y}= \frac{ \log ^2(1+|\xi |^2 )+ \log (1+|\xi |^2 ) + 1 }{1 + \log (1+|\xi |^2 )}\hat{u}, $$
or 
$$ \sqrt{1 + \log (1+|\xi |^2 )}\hat{y}= \frac{ \log ^2(1+|\xi |^2 )+ \log (1+|\xi |^2 ) + 1 }{\sqrt{1 + \log (1+|\xi |^2 )}}\hat{u} . $$
Then 
$$  \| y \| _{Y^1}  =  \int _{{\bf R}_{\xi}^{n}} \left ( 1 + \log (1+|\xi |^2 ) \right )| \hat{y}|^2 d \xi$$
$$ = \int _{{\bf R}_{\xi}^{n}} \frac{\left[ \; 1+  \log ^2(1+|\xi |^2 )+ \log (1+|\xi |^2 )  \; \right] ^2}{1 + \log (1+|\xi |^2 )} | \hat{u} |^2 d \xi.$$
We note that  $$ (1 + \log (1+|\xi |^2 ))^{-1} \left (  \log ^2(1+|\xi |^2 )+ \log (1+|\xi |^2 ) + 1 \right ) ^2 \approx (1+\log(1+|\xi|^2) ^{3}  .$$
Then 
$$ \| y \| _{Y^1} \approx  \int _{{\bf R}_{\xi}^{n}}  (1+\log(1+|\xi|^2) ^{3} | \hat{u} |^2 d \xi = \| u \|_{Y^3}^2. $$
Therefore, $ u \in Y^3 $ and there exists a constant  $c>0$ such that 
$$\| Au \| _{Y^1} \leq c \| u \|_{Y^3} . $$
\hfill $\Box$

\begin{lem} $Y^3 \subset D(A)$.
\end{lem}
{\it Proof.}\
Initially, we observe that $ Y^3  \subset Y^2  $, because $ (1+\log (1+ | \xi |^2) )^2 \leq (1+\log (1+ | \xi |^2) )^3 $. 

Now, let $  u \in Y^3   $, then $ L^{3/2}u \in L^{2}({\bf R}^{n} )$  and we need to  show that there exists  $ y \in Y^1  $ such that
$$ (L^{3/2}u, L^{1/2}\psi ) +(Lu , L \psi ) +(u, \psi) = (y, \psi) + (L^{1/2}y, L^{1/2}\psi),  \text{\; for each \; } \psi \in Y^2 .$$

We define  $ a : Y^1  \times Y^1  \rightarrow {\bf R}  $ by 
$$  a ( \psi, \phi) = ( \psi, \phi) + ( L^{1/2} \psi, L^{1/2} \phi) . $$
This function is well-defined and is a symmetric bilinear form. Moreover, \\
$ \bullet  $ $ a$ is continuous, in fact 
\begin{eqnarray*}
 | a ( \psi , \phi ) |  &\leq & |( \psi, \phi)| +| ( L^{1/2} \psi, L^{1/2} \phi)|  \leq \| \psi \|\vert \phi \| + \| L ^{1/2}\psi \|\vert L^{1/2} \phi \| \\
 &=& \| \hat{ \psi } \|\vert\hat{  \phi }\| + \| \sqrt{\log (1+| \xi | ^2)} \hat{ \psi} \| \| \sqrt{\log (1+| \xi | ^2)} \hat{  \phi } \| \\
 & \leq & 2 \| \sqrt{1+ \log (1+| \xi | ^2)} \hat{ \psi} \| \| \sqrt{1+\log (1+| \xi | ^2)} \hat{ \phi} \| \\
 &=& 2 \| \psi \| _{Y^1} \| \phi \| _{Y^1}.
\end{eqnarray*}
$ \bullet $ $a$ is coercive, because 
$$ a( \phi, \phi ) = ( \phi, \phi) + ( L^{1/2} \phi, L^{1/2} \phi) = \| \phi \|^2 + \| L^{1/2} \phi \|^2 = \|\phi \|_{Y^1}^2 . $$

On the other hand, we define $G: Y^{1}  \rightarrow {\bf R} $ by 
$$ \left \langle G, \psi   \right \rangle = (L^{3/2}u, L^{1/2}\psi ) +(L^{1/2}u , L^{1/2} \psi ) +(u, \psi)  .  $$ 
This map is linear and continuous, in fact 
\begin{eqnarray*}
|\left \langle G, \psi   \right \rangle| & =& |(L^{3/2}u, L^{1/2}\psi )| +|(L^{1/2}u , L^{1/2}\psi| +|(u, \psi)  | \\
&\leq & \|L^{3/2} u \|\|L^{1/2} \psi \| + \|L u \|\| \psi \| + \|u \| \|  \psi \| \\
&=& \|(\log (1+ |\xi |^2))^{3/2} \hat{u} \|\|(\log (1+ |\xi |^2))^{1/2} \hat{\psi } \| + \|\log (1+ |\xi |^2) \hat{u} \| \| \hat{\psi } \| + \|\hat{u} \| \| \hat{\psi} \| \\
&\leq  & \|(\log (1+ |\xi |^2))^{3/2} \hat{u} \| \|(1+ \log (1+ |\xi |^2))^{1/2} \hat{\psi } \| \\ &+& \|\log (1+ |\xi |^2) \hat{u} \| \|(1+ \log (1+ |\xi |^2))^{1/2} \hat{\psi } \| + \|\hat{u} \| \|(1+ \log (1+ |\xi |^2))^{1/2} \hat{\psi } \| \\
&=& \left ( \|(\log (1+ |\xi |^2))^{3/2} \hat{u} \|  +  \|\log (1+ |\xi |^2) \hat{u} \| +  \|\hat{u} \| \right )  \|(1+ \log (1+ |\xi |^2))^{1/2} \hat{\psi } \| \\
&=& \| u \| _{Y^3} \| \psi \| _{Y^{1}}.
\end{eqnarray*}
Thus, we have a continuous linear functional $G$ and a coercive continuous bilinear form in the Hilbert Space $ Y^1 $. By the Lax-Milgram Theorem, there exists a unique  $ y \in Y^{1}  $ such that  
$$ \left \langle G, \psi   \right \rangle = a (y, \psi) \text{\;  for all \; }\psi \in Y^1 . $$
\noindent
In particular, this identity  is valid for all $ \psi \in \mathcal{D}({\bf R}^{n}) $, that is, 
$$ (L^{3/2}u, L^{1/2}\psi ) +(L^{1/2}u , L^{1/2} \psi ) +(u, \psi) = (y, \psi) + (L^{1/2}y, L^{1/2}\psi), \;\;\psi \in \mathcal{D}({\bf R}^{n})   . $$
Since $  \mathcal{D}({\bf R}^{n}) $ is dense in $ Y^2 $, such identity  is valid for all $ \psi \in  Y^2. $ Hence, we can get $ u \in D(A). $
\hfill
$\Box$
\\

Our goal at this section  is to show that $B$ is an infinitesimal generator of a contraction $C^0$-semigroup. For this, we apply the  Lumer-Phillips  theorem. First, we note  that  $D(B)=Y^3 \times Y^2 $ is dense in $ X=Y^2 \times Y^1.  $

 In order to prove that  $B$ is dissipative, we observe that
\begin{eqnarray*}
\left ( B(u,v),(u,v) \right ) _X&=& \left ( (v,-Au), (u,v) \right )_X = (v, u)_{Y^2} + (-Au, v)_{Y^1}  \\
& = &\int _{{\bf R}_{\xi}^{n}} \left ( 1+ \log (1+| \xi |^2) + \log ^2 (1+ |\xi |^2 ) \right ) \hat{v}\overline{\hat{u}} d \xi  - \int _{{\bf R}_{\xi}^{n}}  \left ( 1+ \log (1+| \xi |^2) \right ) \widehat{Au} \overline{\hat{v}} d \xi \\
&=& \int _{{\bf R}_{\xi}^{n}} \left ( 1+ \log (1+| \xi |^2) + \log ^2 (1+ |\xi |^2 ) \right ) \hat{v}\overline{\hat{u}} d \xi  \\ &-& \int _{{\bf R}_{\xi}^{n}}  \left ( 1+ \log (1+| \xi |^2) \right ) \frac{1+ \log (1+| \xi |^2) + \log ^2 (1+ |\xi |^2 )}{1+ \log (1+| \xi |^2)} \hat{u}\overline{\hat{v}} d \xi \\
&=& \int _{{\bf R}_{\xi}^{n}} \left ( 1+ \log (1+| \xi |^2) + \log ^2 (1+ |\xi |^2 ) \right ) (  \hat{v}\overline{\hat{u}}  - \hat{u}\overline{\hat{v}}) d \xi \\
&=& 2i \int _{{\bf R}_{\xi}^{n}} \left ( 1+ \log (1+| \xi |^2) + \log ^2 (1+ |\xi |^2 ) \right )  \text{Im} (\hat{u}\overline{\hat{v}}) d \xi.  
\end{eqnarray*}
Thus, $ \text{Re} \left ( B(u,v),(u,v) \right ) _X = 0 $ and $B$ is dissipative. 
\vspace{0.2cm}

Now we need to prove that $(I - B)(X) = X$. Clearly $ (I - B)D(B)  \subset X  $. In its turn, let $ (f, g ) \in X = Y^2 \times Y^1  $. Then, we will prove that there exists $ (u,v) \in D(B) $ such that $ (I-B)(u,v) = (f,g) .$ \\
We define a mapping ${\bf  a} : Y^2 \times Y^2 \rightarrow {\bf R}$ by
$$ {\bf a}(\varphi, \psi ) = (\varphi , \psi )_{Y^1} + (\varphi , \psi )_{Y^2} .$$
Then $ {\bf a}  $ is a symmetrical bilinear form,  which is \\
$ \bullet $ continuous:
$$ | {\bf a}(\varphi, \psi ) | \leq \| \varphi \|_{{Y^1}} \| \psi \|_{{Y^1}} + \| \varphi \|_{{Y^2}} \| \psi \|_{{Y^2}} \leq 2  \| \varphi \|_{{Y^2}} \| \psi \|_{{Y^2}},$$
$ \bullet $ coercive: 
$$ {\bf a}(\psi, \psi ) = (\psi , \psi )_{Y^1} + (\psi , \psi )_{Y^2} = \| \psi \|_{Y^1} ^2+ \| \psi \|_{Y^2} ^2 \geq \| \psi \|_{Y^2} ^2 .$$
\noindent
Furthermore, we consider the linear functional $ F:  Y^2 \rightarrow {\bf R}$ given by 
$$ \left \langle  F, \psi \right \rangle = (f+g, \psi )_{Y^1}, $$
which is well-defined because of $ Y^2 \subset Y^1  . $ Also, one has 
$$ | \left \langle  F, \psi \right \rangle | \leq \| f+g \|_{Y^1} \| \psi \|_{Y^1} \leq  \| f+g \|_{Y^1} \| \psi \|_{Y^2} ,$$
which just proves the continuity of $ F$. 

Thus, we can apply the  Lax-Milgram theorem to get the existence of a unique $ u \in Y^2 $ such that
$$ {\bf a}(u, \psi ) = F( \psi ) \text{ for all } \psi \in Y^2 . $$
In other words, 
$$ (u, \psi ) + (L^{1/2}u, L^{1/2}\psi ) + (u, \psi ) + (L^{1/2}u, L^{1/2}\psi ) + (Lu, L \psi )$$
$$ = (f+g, \psi ) +(L^{1/2}(f+g), L^{1/2}\psi ) $$
 for all $ \psi \in Y^2  $. In particular, this equality is valid for all $ \psi \in \mathcal{D}({\bf R}^{n})  $  and we have the following  identity  in $ \mathcal{D}'({\bf R}^{n}) $
 $$ A u+ u = f+g . $$

Finally, we observe that $ u \in Y^3 . $ In fact, by applying the Fourier transform, we obtain $$ \widehat{Au} + \hat{u} = \hat{f} + \hat{g}$$
and 
 $$ 	\frac{1 + \log (1+|\xi |^2 ) + \log ^2 (1+|\xi |^2 )}{\log (1+|\xi |^2 )} \hat{u} + \hat{u} = \hat{f}+\hat{g},$$
which is equivalent to 
 $$ \frac{1 + \log (1+|\xi |^2 ) + \log ^2 (1+|\xi |^2 )}{ \sqrt{1+ \log (1+|\xi |^2 )}} \hat{u}  = \sqrt{1+ \log (1+|\xi |^2 )} (\hat{f}+\hat{g} - \hat{u}). $$
 
From the fact that
$$ \frac{1 + \log (1+|\xi |^2 ) + \log ^2 (1+|\xi |^2 )}{ \sqrt{1+ \log (1+|\xi |^2 )}}  \approx (1+\log (1+|\xi |^2) )^{3/2},$$
we have
$$ (1+\log(1+|\xi |^2))^3 | \hat{u} |^2 \approx (1+\log(1+|\xi |^2)) | \hat{f} + \hat{g} - \hat{u}|^2.$$
Now, since $ f, g, u \in Y^{1} $ we can conclude
 $$  \int _{{\bf R}_{\xi}^{n}} (1+\log(1+|\xi |^2))^3 | \hat{u}(\xi) |^2 d \xi  < \infty .$$
This estimate proves that $u \in Y^3$. Then,  $ v = u - f \in Y^2  $ because of $f \in Y^2$ and we have obtained that
 $$ (I-B) (u,v) = (f,g) .$$
 Hence, we have proved  that  $(I-B)D(B) = X.$ Therefore, by the Lumer-Phillips theorem the operator $B$ is an infinitesimal generator of a $C^0-$semigroup of contractions.
 \vspace{0.1cm}
 
 Next we want to prove that $  J: X \rightarrow X $ is  a linear bounded operator. The linearity is simple. The boundedness is given by the following series of inequalities:
\begin{eqnarray*}
\| J(u,v) \|^2_X 
 &=& \left \| \sqrt{1+\log(1+ |\xi |^2 )} \frac{\hat{u }- \hat{v}}{ 1+\log(1+ |\xi |^2 ) } \right \|^2 \\
&\leq & 2 \int _{{\bf R}_{\xi}^{n}}\frac{| \hat{ u}|^2 }{ 1+\log(1+ |\xi |^2 ) } d \xi + 2 \int _{{\bf R}_{\xi}^{n}}\frac{| \hat{ v}|^2 }{ 1+\log(1+ |\xi |^2 ) } d \xi \\
&\leq & 2\int _{{\bf R}_{\xi}^{n}} | \hat{ u}|^2 d \xi + 2 \int _{{\bf R}_{\xi}^{n}}  ( 1+\log(1+ |\xi |^2 ) ) | \hat{ v}|^2 d \xi \\
&\leq & 2 \|u \|_{Y^2}^2 + 2 \| v \| _{Y^1}^{2}\\
&=& 2 \| (u,v) \|_{X}^2.
\end{eqnarray*}
Therefore, $ B+J $ is an infinitesimal generator of a $ C^0-$semigroup in $X$ and we have arrived at the following result.
\begin{theo}\label{existence}
Let $ n \geq 1 $ and  $ (u _0, u_1)  \in Y^3 \times Y^2$. Then the problem \eqref{eqn}--\eqref{initial} admits a unique  solution in the class
$$ u \in C([0, \infty), Y^3) \cap C^1([0, \infty), Y^2) \cap C^2([0, \infty), Y^1). $$
Moreover, for initial data  $(u_0,u_1) \in X=Y^2 \times Y^1$  the problem \eqref{eqn}--\eqref{initial} admits a unique weak solution in the class
$$ u \in C([0, \infty), Y^2) \cap C^1([0, \infty), Y^1) \cap C^2([0, \infty), L^2). $$
\end{theo}
\section{Asymptotic behavior via multiplier method}

We begin with this section by considering the  Cauchy problem in the Fourier space associated to the problem \eqref{eqn}--\eqref{initial}  as follows 
\begin{align}\label{Eq5}
& (1+\log (1+|\xi|^2)) \hat{u}_{tt} + \log (1+|\xi|^2) (1+\log (1+|\xi|^2)) \hat{u} + \hat{u}_t = 0, \nonumber \\
& \hat{u}(0, \xi) =  \hat{u}_0(\xi),  \\ 
& \hat{u}_t(0, \xi) =  \hat{u}_1(\xi).  \nonumber 
\end{align}

Multiplying the equation in \eqref{Eq5} by $\overline{\hat{u}_t}  $ one can get the following pointwise energy identity
\begin{equation}\label{Eq8}
\frac{\mathrm{d} }{\mathrm{d} t} E_0(t, \xi) + | \hat{u}_t(t,\xi)|^2=0,
\end{equation}
where 
\begin{equation*}
E_0(t, \xi) =\frac{(1+\log(1+|\xi|^2)) | \hat{u}_t |^2 }{2} + \frac{\log (1+|\xi|^2)(1+\log(1+|\xi|^2))| \hat{u}|^2}{2}, \quad t>0, \; \xi \in {\bf R}^{n} 
\end{equation*}
 is the total density of the energy for the system  \eqref{Eq5}.  From \eqref{Eq8} we see that $E_0(t, \xi)$ is a  non-increasing function of  $t$ for each $\xi$.
 
\begin{lem}\label{lemmadefinicaorho} Consider the following functions  
$$ \varphi(\xi ) = \log (1+|\xi|^2)(1+\log(1+|\xi|^2)), \quad \psi(\xi)= \frac{1}{1+\log(1+|\xi|^2)},  \quad\phi(\xi) = \sqrt{\log(1+|\xi|^2)}.  $$
Then, there is a unique real number $ \delta _0 \in (0,1)$ such that  
\begin{itemize}
\item[{\rm (i)}.]$|\xi | \leq \delta _0 $ $\Rightarrow$ $\varphi (\xi) \leq \psi (\xi) $ and $ \varphi (\xi) \leq \phi (\xi) $,
\item[{\rm(ii)}.] $|\xi | \geq \delta _0 $ $\Rightarrow$ $\psi (\xi ) \leq \varphi (\xi) $ and  $ \psi (\xi ) \leq \phi (\xi) $.
\end{itemize}
\end{lem}
{\it Proof.}\
Let   $ \theta (r ) = (1+\log(1+r^2)) \sqrt{\log (1+r^2) } -1 $ for $ r \geq 0 $. We note  that $\theta(0)= -1$ and $ \theta(1) = (1+\log 2) \sqrt{\log 2} -1 >0  $.  The continuity of this function implies there exists $ 0<\delta _0<1 $ such that $ \theta (\delta _0)=0 $. Moreover, $\theta (r)$ is a increasing function of  $r\geq 0$ and such fact implies that  $\delta _0 >0$ is a unique number satisfying $-1 \leq  \theta (r) \leq 0  $ for all  $ 0 \leq r \leq \delta _0 $ and $ \theta (r) >0 $ for  $ r >\delta _0 $.
\noindent
That is,  if $ 0  \leq  r  \leq \delta  _0 $, then one has 
$$ (1+\log(1+r^2)) \sqrt{\log (1+r^2) } \leq 1 . $$
Thus, $$  (1+\log(1+r^2)) \log (1+r^2)  \leq \sqrt{\log (1+r^2) } \text{ and }  \sqrt{\log (1+r^2) } \leq \frac{1}{1+\log(1+r^2)}\quad (0 \leq r \leq \delta_{0}). $$
Similarly, if $ r > \delta _0$,  then
$$ (1+\log(1+r^2)) \sqrt{\log (1+r^2) } > 1,$$
and one obtains 
$$ (1+\log(1+r^2)) \log (1+r^2)  > \sqrt{\log (1+r^2) } \text{ and }  \sqrt{\log (1+r^2) } > \frac{1}{1+\log(1+r^2)}\quad (r > \delta_{0}).$$
These imply the desired estimates (i) and (ii).   
\hfill
$\Box$
\\

For $\delta _0 >0 $ given in Lemma \ref{lemmadefinicaorho}, we define the following function of $ \xi \in {\bf R}^{n} $ such that
\begin{equation}\label{funcaorho}
\rho (\xi) := \left\{\begin{matrix}
\frac{1}{2} \log (1+|\xi|^2)(1+\log(1+|\xi|^2)), & \text{ for } | \xi | \leq \delta _0, \\ 
\hspace{-2.5cm} \frac{1}{2(1+\log(1+|\xi|^2))}, & \text{ for } | \xi | > \delta _0 .
\end{matrix}\right.
\end{equation}
As a consequence of Lemma \ref{lemmadefinicaorho}, we have
$$ \rho (\xi ) = \text{min}  \left \{ \frac{1}{2} \log (1+|\xi|^2)(1+\log(1+|\xi|^2)), \frac{1}{2(1+\log(1+|\xi|^2))}  , \frac{1}{2}  \sqrt{\log(1+|\xi|^2)}   \right \} .$$
By multiplying the equation \eqref{Eq5} by $ \rho (\xi) \overline{\hat{u}} $ we obtain the identity 
\begin{equation}\label{Eq10}
\rho (\xi)(1+\log (1+|\xi|^2)) \hat{u}_{tt}  \overline{\hat{u}} + \rho (\xi)\log (1+|\xi|^2) (1+\log (1+|\xi|^2)) |\hat{u} |^2 + \frac{\rho(\xi)}{2} \frac{\mathrm{d} }{\mathrm{d} t} |\hat{u}|^2 = 0
\end{equation}
for all $t >0 $ and $ \xi \in {\bf R}^{n}  $.  Taking the real part on the last identity we arrive at
\begin{align}\label{Eq11}
& \frac{\mathrm{d} }{\mathrm{d} t} \left [  \rho (\xi)(1+\log (1+|\xi|^2)) \text{Re}( \hat{u}_{t}  \overline{\hat{u}} ) +  \frac{ \rho(\xi) }{2} |\hat{u}|^2 \right ] + \rho(\xi) \log(1+|\xi|^2 ) (1+\log (1+|\xi|^2)) |\hat{u} |^2  \nonumber \\ &= \rho (\xi) (1+\log (1+|\xi|^2)) |\hat{u}_t |^2  .
\end{align}
To proceed further we need to define the following functions on $ (0, \infty ) \times  {\bf R}_{\xi}^{n} $ such that 
\begin{align}
& E(t, \xi) = E_0(t, \xi) + \rho (\xi) (1+\log (1+|\xi|^2))\text{Re}(\hat{u}_t \overline{\hat{u}} )+ \frac{\rho (\xi)}{2} | \hat{u} |^2, \nonumber \\
& F(t, \xi) = |\hat{u}_t|^2 + \rho (\xi )\log (1+|\xi | ^2)(1+ \log (1+|\xi | ^2))| \hat{u}|^2, \\
& R(t, \xi) = \rho (\xi) (1+\log (1+|\xi|^2)) |\hat{u}_t|^2. \nonumber
\end{align}
Then, by adding \eqref{Eq8} and \eqref{Eq11}, we get the following identity
\begin{equation}\label{3.7}
\frac{\mathrm{d} }{\mathrm{d} t} E(t, \xi) + F(t, \xi) = R(t,\xi), \quad t>0,\; \xi \in {\bf R}_{\xi}^{n}.
\end{equation}

\begin{lem}\label{lemmaequivalenciaenergia}\, It holds that
 $$ \frac{1}{2} E_0(t, \xi) \leq E(t, \xi ) \leq 3 E_0(t, \xi), \text{ } t >0, \text{ } \xi \in {\bf R}_{\xi}^{n}.  $$
\end{lem}
{\it Proof.}\ By Lemma \ref{lemmadefinicaorho} we have for $t>0$ and $\xi \in {\bf R}_{\xi}^{n}$,
\begin{eqnarray*}
E(t , \xi )
&\leq & E_0(t, \xi ) + \rho (\xi )(1+ \log (1+|\xi |^2) )|\hat{u}_t | |\hat{u}|+ \frac{\rho (\xi )}{2}| \hat{u}|^2 \\
&\leq & E_0(t, \xi) + \frac{1+ \log (1+|\xi |^2)}{2}|\hat{u}_t |^2 + \frac{\rho (\xi)^2 (1+ \log (1+|\xi |^2))}{2} | \hat{u}|^2 + \frac{\rho (\xi )}{2}| \hat{u}|^2 \\
&\leq & E_0(t, \xi) + \frac{1+ \log (1+|\xi |^2)}{2}|\hat{u}_t |^2 +  \frac{ \log(1+|\xi | ^2) (1+ \log (1+|\xi |^2))}{8} | \hat{u}|^2 \\ 
&+&  \frac{ \log(1+|\xi | ^2) (1+ \log (1+|\xi |^2)) }{4}| \hat{u}|^2 \\ 
&\leq & 3 E_0(t, \xi),
\end{eqnarray*}
according to the definition of $E_0(t,\xi)$ in \eqref{Eq8}.

On the other hand, one has
\begin{eqnarray*}
- \rho (\xi ) (1+\log (1+|\xi| ^2)) \text{Re} (\hat{u}_t\overline{\hat{u}}) &\leq & \rho (\xi ) (1+\log (1+|\xi| ^2))  |\hat{u} _t| |\hat{u}| \\
&\leq & \frac{1+\log (1+|\xi| ^2)}{4}|\hat{u} _t|^2 + \rho(\xi)^2(1+\log (1+|\xi| ^2))| \hat{u}|^2 \\
&\leq & \frac{1+\log (1+|\xi| ^2)}{4}|\hat{u} _t|^2 + \frac{1}{4} \log(1+|\xi |^2)(1+\log (1+|\xi| ^2))| \hat{u}|^2. 
\end{eqnarray*}
\noindent
Thus, the last estimate implies  
\begin{eqnarray*}
E(t, \xi) &=& E_0(t, \xi ) + \rho (\xi )(1+ \log (1+|\xi |^2) )\text{Re} (\hat{u}_t \overline{\hat{u}}) + \frac{\rho (\xi )}{2}| \hat{u}|^2 \\
&\geq & E_0(t, \xi ) + \rho (\xi )(1+ \log (1+|\xi |^2) )\text{Re} (\hat{u}_t \overline{\hat{u}}) \\
&\geq & \left ( \frac{1}{2} - \frac{1}{4} \right )(1+ \log (1+|\xi |^2) ) | \hat{u}_t |^2+ \left ( \frac{1}{2} - \frac{1}{4} \right )\log(1+|\xi |^2)(1+\log (1+|\xi| ^2))| \hat{u}|^2 \\
&=&  \frac{1}{4} (1+ \log (1+|\xi |^2) ) | \hat{u}_t |^2+  \frac{1}{4} \log(1+|\xi |^2)(1+\log (1+|\xi| ^2))| \hat{u}|^2 \\
&=& \frac{1}{2} E_0(t, \xi),
\end{eqnarray*}
for $t>0$ and $\xi \in {\bf R}_{\xi}^{n}$. These imply the desired estimates.
\hfill
$\Box$

Now, we need the next lemma.
\begin{lem}\label{lemma'gronwall'}\, It is true that
$$ \frac{\mathrm{d} }{\mathrm{d} t} E(t, \xi) + \frac{\rho(\xi )}{2} E(t, \xi ) \leq 0 , \text{ } t >0, \text{ } \xi \in {\bf R}_{\xi}^{n}.$$
\end{lem}
{\it Proof.}\ By definitions of $E(t,\xi)$, $ F(t,\xi) $, $R(t,\xi)$, $\rho(\xi)$ and \eqref{3.7}, one obtains a series of inequalities 
\begin{eqnarray*}
\frac{\mathrm{d} }{\mathrm{d} t} E(t, \xi) + \frac{\rho(\xi )}{2} E(t, \xi )  &=& R(t,\xi) -  F(t, \xi) +  \frac{\rho(\xi )}{2} E(t, \xi ) \\
&\leq &  R(t,\xi) -  F(t, \xi) +  \frac{3 \rho(\xi )}{2} E_0(t, \xi ) \\
&=&  \left ( \frac{7}{4} \rho (\xi) (1+\log(1+|\xi|^2)) -1  \right )|\hat{u}_t|^2 \\ &-& 
\frac{1}{4} \rho (\xi) \log (1+|\xi|^2)(1+\log(1+|\xi|^2)) | \hat{u}|^2 \\
&\leq &  - \frac{1}{8} |\hat{u}_t|^2 - \frac{1}{4} \rho (\xi) \log (1+|\xi|^2)(1+\log(1+|\xi|^2)) | \hat{u}|^2 \\
&\leq & 0
\end{eqnarray*}
for $t>0$ and ${\bf R}_{\xi}^{n}$. This implies the desired statement.
\hfill
$\Box$

Now, it follows from Lemma \ref{lemma'gronwall'} that
$$ E(t, \xi ) \leq E(0, \xi ) e^{-\frac{\rho (\xi)}{2}t} ,$$
for $t>0$ and $\xi \in {\bf R}_{\xi}^{n}$.
\noindent
By combining the last estimate with Lemma \ref{lemmaequivalenciaenergia}, one can deduce the important pointwise estimate,
$$ E_0(t, \xi ) \leq 6 E_0(0, \xi ) e^{-\frac{\rho (\xi)}{2}t} , $$
for $t>0$ and $\xi \in {\bf R}_{\xi}^{n}$.

The above estimate combined with the definition of $E_0(t,\xi)$ in \eqref{Eq8} implies the following crucial pointwise estimate.
\begin{pro}\label{propestimativaUmetododeenergia} It holds that 
\begin{align*}
& (1+\log(1+|\xi|^2)) | \hat{u}_t(t,\xi)|^2 + \log (1+|\xi|^2)(1+\log(1+|\xi|^2))| \hat{u}(t,\xi)|^2 \\ &\leq 6 (1+\log(1+|\xi|^2)) e^{-\frac{\rho (\xi)}{2}t} | \hat{u}_1(\xi)|^2 + 6 \log (1+|\xi|^2)(1+\log(1+|\xi|^2)) e^{-\frac{\rho (\xi)}{2}t} | \hat{u}_0(\xi)|^2 
\end{align*}
for all $t >0 $ and $ \xi \in {\bf R}_{\xi}^{n}$, and
$$ | \hat{u}(t,\xi)|^2 \leq 6 \frac{ e^{-\frac{\rho (\xi)}{2}t}  }{\log (1+|\xi |^2)} | \hat{u}_1(\xi)|^2 + 6 e^{-\frac{\rho (\xi)}{2}t}  | \hat{u}_0(\xi)|^2 $$ for all $t >0 $ and $ \xi \in {\bf R}_{\xi}^{n}$, $ \xi \neq 0 $. 
\end{pro}


As a consequence of the second estimate in Proposition \ref{propestimativaUmetododeenergia} one can get the following result.
\begin{pro}\,Let $n \geq 3$, and let $u(t,\xi)$ be the solution to problem \eqref{eqn}-\eqref{initial}.  Suppose $u_0 \in L^1({\bf R}^{n}) \cap Y^{\frac{n-2}{2}}$, $u_1 \in L^1({\bf R}^{n}) \cap Y^{\frac{n-2}{2}}$. Then,  the following estimate holds:
$$ \int _{{\bf R}_{\xi}^{n}} | u(t, x) |^2 d x   \leq   C_1 t^{-\frac{n-2}{2}} \Big[ \|  u_1 \|_{L^{1}}^2 + 
 \| u_0 \| _{Y^{\frac{n-2}{2}}} ^2\Big]
+ 
C_2  t^{- \frac{n}{2}} \Big[\| u_1 \| _{Y^{\frac{n-2}{2}}} ^2   +  \|  u_0 \|_{L^{1}}^2  \Big] , \;\; t>0,$$
where $ C_1$ and $C_2 $ are positive constants depending only on $n$. 
\end{pro}
{\it Proof.}\  Let $\delta_0 > 0$ be a given real number obtained in Lemma \ref{lemmadefinicaorho} . To prove the Proposition above one needs to consider separately the zones of low and high frequency.

\textbf{1) Estimate on the zone $|\xi | \leq \delta _0 $}

On this region one notices $\rho (\xi) = \frac{1}{2} \log (1+|\xi |^2 ) (1+  \log (1+|\xi |^2))$. Then, one can observe that $1 \leq 1+  \log (1+|\xi |^2) \leq  1+  \log (1+\delta _0 ^2)$ for $|\xi|\leq \delta_0$. Thus we get 
$$\frac{1}{2}\log(1+\vert\xi\vert^{2}) \leq \rho (\xi) \leq \frac{1+  \log (1+\delta  _0 ^2)}{2}\log (1+\delta  _0 ^2),\quad| \xi | \leq \delta  _0.$$ 
Then, by applying the second estimate of Proposition \ref{propestimativaUmetododeenergia}, one obtains 
\begin{eqnarray*}
\int _{|\xi | \leq \delta  _0} | \hat{u}|^2 d \xi &\leq& 6 \int _{|\xi | \leq \delta  _0} \frac{ e^{-\frac{\rho (\xi)}{2}t}  }{\log (1+|\xi |^2)} | \hat{u}_1 |^2 d \xi + 6 \int _{|\xi | \leq \delta  _0} e^{-\frac{\rho (\xi)}{2}t}  | \hat{u}_0|^2 d \xi \\
&\leq & 6 \int _{|\xi | \leq \delta  _0} \frac{ e^{-\frac{\log(1+|\xi|^2)}{4}t}  }{\log (1+|\xi |^2)} | \hat{u}_1 |^2 d \xi + 6 \int _{|\xi | \leq \delta  _0} e^{-\frac{\log(1+|\xi|^2)}{4}t}  | \hat{u}_0|^2 d \xi \\
&= & 6 \int _{|\xi | \leq \delta  _0} (1+|\xi|^2)^{-\frac{t}{4}}\frac{ 1 }{\log (1+|\xi |^2)} | \hat{u}_1 |^2 d \xi + 6 \int _{|\xi | \leq \delta  _0}   (1+|\xi|^2)^{-\frac{t}{4}} | \hat{u}_0|^2 d \xi \\
&\leq & 6  \|  \hat{u}_1 \|_{L^{\infty}}^2 \int _{|\xi | \leq \delta  _0} (1+|\xi|^2)^{-\frac{t}{4}} \frac{ 1 }{\log (1+|\xi |^2)} d \xi + 6 \| \hat{u}_0 \|_{L^{\infty}}^2 \int _{|\xi | \leq \delta  _0}   (1+|\xi|^2)^{-\frac{t}{4}}  d \xi \\
&\leq & 6 \omega _n  \|  u_1 \|_{L^{1}}^2 \int _{0}^{\delta  _0} (1+r^2)^{-\frac{t}{4}} \frac{ r^2 }{\log (1+r^2)} r^{n-3}d r + 6 \omega _n \| u_0 \|_{L^{1}}^2 \int _{0}^{\delta  _0}  (1+r^2)^{-\frac{t}{4}} r^{n-1}  d r \\
&\leq & C_{0,n} \omega _n  \|  u_1 \|_{L^{1}}^2 \int _{0}^{\delta  _0} (1+r^2)^{-\frac{t}{4}}  r^{n-3}d r + 6 \omega _n \| u_0 \|_{L^{1}}^2 \int _{0}^{\delta  _0}  (1+r^2)^{-\frac{t}{4}} r^{n-1}  d r \\
&\leq & C_{0,n}\omega _n  \|  u_1 \|_{L^{1}}^2 \int _{0}^{1} (1+r^2)^{-\frac{t}{4}}  r^{n-3}d r + 6 \omega _n \| u_0 \|_{L^{1}}^2 \int _{0}^{1}  (1+r^2)^{-\frac{t}{4}} r^{n-1}  d r \\
&\leq & C_{1,n}  \|  u_1 \|_{L^{1}}^2 t^{-\frac{n-2}{2}}  +  C_{2,n}  \|  u_0 \|_{L^{1}}^2 t^{-\frac{n}{2}}, \;\; t>0
\end{eqnarray*}
with some constants $C_{j,n} > 0$ ($j = 0,1$), where we have just used  Lemma \ref{general-p} and the fact that
\[\lim_{\sigma \to +0}\frac{\sigma}{\log(1+\sigma)} = 1.\]
\\

\textbf{2) Estimate on the zone $|\xi | \geq \delta_0$} \\ 
In this case, one notices $\rho (\xi ) =  \displaystyle{\frac{1}{2 (1+\log (1+|\xi|^2))}}$. By  Proposition \ref{propestimativaUmetododeenergia} and the definition of $\rho(\xi)$ we have

\begin{eqnarray*}
\int _{|\xi | [geq \delta  _0} | \hat{u}|^2 d\xi &\leq & 6 \int _{|\xi | \geq \delta  _0} \frac{ e^{-\frac{\rho (\xi)}{2}t}  }{\log (1+|\xi |^2)} | \hat{u}_1 |^2 d \xi + 6 \int _{|\xi | \geq \delta  _0} e^{- \frac{\rho (\xi)}{2}t}  | \hat{u}_0|^2d \xi  \\
&=&    6 \int _{|\xi | \geq \delta  _0} \frac{ e^{-\frac{1}{4(1+\log(1+|\xi|^2))} t}  }{\log (1+|\xi |^2)} | \hat{u}_1 |^2 d \xi + 6 \int _{|\xi | \geq \delta  _0} e^{ -\frac{1}{4(1+\log(1+|\xi|^2))} t}  | \hat{u}_0|^2d \xi. 
\end{eqnarray*}
Now, applying Lemma \ref{LemmaTecnico1} with $\nu =\frac{n}{2} $ and $a = -1$, and  $\nu' =\frac{n-2}{2}$ and $a = -1$ to the last two integrals above, one has 
\begin{eqnarray*} 
\int _{|\xi | > \delta  _0} | \hat{u}|^2 d\xi 
&\leq & C t^{-\nu} \int _{|\xi | \geq \delta _0}\hspace{-0.2cm} \frac{ (1+\log(1+| \xi |^2))^{\nu}  }{\log (1+|\xi |^2)} | \hat{u}_1 |^2 d \xi + C t^{-\nu'} \int _{|\xi | \geq \delta  _0} \hspace{-0.1cm}(1+\log(1+| \xi |^2))^{\nu'}  | \hat{u}_0|^2d \xi \\
&\leq & C t^{- \frac{n}{2}} \int _{|\xi | \geq \delta  _0} \frac{ (1+\log(1+| \xi |^2))^{\frac{n}{2}}  }{\log (1+|\xi |^2)} | \hat{u}_1 |^2 d \xi \\ &+& C t^{- \frac{n-2}{2}} \int _{|\xi |\geq \delta  _0} (1+\log(1+| \xi |^2))^{\frac{n-2}{2} }  | \hat{u}_0|^2d \xi \\
&=& C t^{- \frac{n}{2}} \int _{|\xi | \geq \delta  _0} \frac{  1+ \log (1+|\xi |^2) }{\log (1+|\xi |^2)}  (1+\log(1+| \xi |^2))^{\frac{n-2}{2}} | \hat{u}_1 |^2 d \xi \\ &+& C t^{- \frac{n-2}{2}} \int _{|\xi | \geq \delta  _0} (1+\log(1+| \xi |^2))^{\frac{n-2}{2} }  | \hat{u}_0|^2d \xi \\
&\leq & C_1 t^{- \frac{n}{2}} \int _{|\xi | \geq \delta  _0}  (1+\log(1+| \xi |^2))^{\frac{n-2}{2}} | \hat{u}_1 |^2 d \xi \\
&+& C t^{- \frac{n-2}{2}} \int _{|\xi | \geq \delta  _0} (1+\log(1+| \xi |^2))^{\frac{n-2}{2} }  | \hat{u}_0|^2d \xi \\
&\leq & C_1 t^{- \frac{n}{2}} \| u_1 \| _{Y^{\frac{n-2}{2}}}^2 + C t^{- \frac{n-2}{2}}  \| u_0 \| _{Y^{\frac{n-2}{2}}}^2,
\end{eqnarray*}
where one has just used the property 
\[\lim_{\sigma \to \infty}\frac{1+\log(1+\sigma)}{\log(1+\sigma)} = 1.\]
By adding two estimates for low and high frequencies and using the Plancherel theorem, one can conclude the proof of proposition.  
\hfill
$\Box$


\section{Asymptotic profile of solutions}

Applying the Fourier transform on the problem \eqref{eq2}  one obtain the associated problem in Fourier space:
\begin{align}\label{Probl-Fourier}
& (1+\log (1+ |\xi|^2)) \hat{u}_{tt} + \log (1+|\xi|^2) (1+\log (1+|\xi|^2)) \hat{u} + \hat{u}_t = 0, \nonumber\\
& \hat{u}(0, \xi) = \hat{u}_0(\xi),\\ 
& \hat{u}_t(0, \xi) = \hat{u}_1(\xi).\nonumber
\end{align}
 The  characteristic roots of the associated  polynomial  to the  equation in \eqref{Probl-Fourier} are given by
$$ \lambda _{\pm } =  \frac{-1 \pm \sqrt{1-4 \log (1+|\xi|^2 )(1+\log (1+|\xi|^2))^2}}{2 (1+ \log (1+|\xi|^2))}.$$

 We observe that there exists a unique real number $ \delta >0 $ such that 
\begin{equation}\label{condicao}
 1-4 \log (1+|\xi|^2 )(1+\log (1+|\xi|^2))^2 \left\{\begin{matrix}
\geq 0 & \text{ for} |\xi| \leq \delta, \\ 
< 0  &  \text{ for} | \xi | > \delta. 
\end{matrix}\right. 
\end{equation}
In fact, the function $ f(r) = 1-4 \log (1+r^2 )(1+\log (1+r^2))^2 $ is decreasing for $ r \geq 0$, continuous and 
$$f (0) = 1, f(1) = 1 - 4 \log 2\cdot(1+\log 2)^2 <0 .$$
Therefore,  by the mean value theorem there exists a unique number $\delta  \in (0,1)$ that satisfies \eqref{condicao}. The same theorem guarantees us the existence of  $ 0< \eta < \delta  $ such that 
\begin{equation}\label{def-eta}
 \frac{1}{2}\leq \sqrt{1 - 4 \log (1+ |\xi |^2 ) (1+\log (1+|\xi |^2))^2} \leq 1, 
 \end{equation}
whenever  $ | \xi | \leq \eta .$ 

The  next lemma is very important to  get sharp estimates. In the following notation $A \approx B$ means that  $c_1 A \leq B \leq c_2 A$ for some positive constants $c_1, c_2$.
\begin{lem}\label{lemadeequivalencias} It holds that
\begin{itemize}
\item[{\rm (i)}.] $ \lambda _+ \approx - \log (1+|\xi|^2)$,
\item[{\rm (ii)}.]$\lambda _- \approx -1 $,
\item[{\rm (iii)}.] $ \lambda _+ +\lambda _- \approx -1 $,
\end{itemize}
whenever $|\xi | \leq \delta$. And, in particular, in the case of $ |\xi | \leq \eta $, one has
\begin{itemize}
\item[{\rm (iv)}.]$\lambda _+ - \lambda _- \approx 1 .$
\end{itemize}
\end{lem} 
\begin{rem}
{\rm Note that the constants $c_1$ and $c_2$ appeared in Lemma \ref{lemadeequivalencias} may depend on $\delta$ or $\eta$. We also note that the four items in Lemma \ref{lemadeequivalencias} simultaneously hold on the zone $\{|\xi| \leq \eta\}$ because of $\eta<\delta$.}
\end{rem}
{\it Proof.}\
\begin{itemize}
\item[{\rm (i)}.]\,First, we observe that  $ 0 \leq 4 \log ^2 (1+|\xi |^2 ) (1+ \log (1+|\xi |^2 ))^4 $ for all $ \xi \in {\bf R}_{\xi}^{n}$. For this reason, one has
\begin{align*}
1 - 4 \log (1+|\xi |^2 ) (1+ \log (1+|\xi |^2 ))^2 &\leq  \left [ 1- 2 \log (1+|\xi |^2 )(1+ \log (1+|\xi |^2 ))^2 \right ]^2. 
\end{align*}
And, for $ |\xi | \leq \delta$ we have 
\begin{align*}
\sqrt{1 - 4 \log (1+|\xi |^2 ) (1+ \log (1+|\xi |^2 ))^2} \leq 1- 2 \log (1+|\xi |^2 )(1+ \log (1+|\xi |^2 ))^2.
\end{align*}
This estimate implies that $$  \lambda _+ = \frac{-1 + \sqrt{1- 4 \log (1+|\xi|^2)(1+ \log (1+|\xi|^2))^2} }{2 (1+ \log(1+|\xi |^2))} \leq - \log (1+|\xi|^2)(1+\log(1+|\xi|^2)) .  $$
Now, because of $ 1 \leq 1+\log(1+|\xi|^2) $ we obtain the upper bound to $\lambda_+$ in the case of $\vert\xi\vert \leq \delta$:
 $$ \lambda _+ \leq - \log (1+|\xi |^2).$$
 
 On the other hand,  for $ |\xi | \leq \delta $, one has $ 0\leq 1- 4 \log (1+|\xi | ^2) (1+\log (1+|\xi |^2))^2 \leq 1$, so that 
 \begin{align*}
  1 - 4 \log (1+|\xi |^2 ) (1+ \log (1+|\xi |^2 ))^2 \leq \sqrt{1 - 4 \log (1+|\xi |^2 ) (1+ \log (1+|\xi |^2 ))^2}. 
 \end{align*}
Hence one can get the following lower bound on $\lambda_{+}$ such that
$$  -2\log (1+|\xi | ^2) (1+\log (1+|\xi |^2))  \leq \frac{-1 + \sqrt{1 -4 \log (1+|\xi | ^2) (1+\log (1+|\xi |^2))^2} }{ 2 (1+\log (1+|\xi |^2))} = \lambda _+ .$$
Therefore, we can get 
$$ -2K_{\delta} \log (1+|\xi | ^2)   \leq  \lambda _+ , \quad |\xi|\leq \delta,$$
with $K_{\delta}=: 1+\log(1+\delta^2)$.
\noindent
These imply the desired statement of {\rm (i)}.

  \item[{\rm (ii)}.] Since  $0 \leq \sqrt{1-4 \log (1+|\xi|^2 )(1+\log (1+|\xi|^2))^2} \leq 1$ in the region $|\xi|\leq \delta$, 
$$  \frac{-1}{1+ \log (1+|\xi |^2)} \leq \frac{-1 - \sqrt{1-4 \log (1+|\xi|^2 )(1+\log (1+|\xi|^2))^2}}{2(1+\log (1+|\xi |^2))} \leq \frac{-1}{2(1+\log (1+|\xi |^2))} .$$  
Therefore,  $$ -1\leq \lambda _- \leq \frac{-1}{2K_{\delta}}.$$
 
  \item[{\rm (iii)}.]  To prove this item we observe that $ \lambda _+ + \lambda _- = \frac{-1}{1+ \log (1+|\xi |^2)} $. Hence, 
  $$ -1 \leq \lambda _+ + \lambda _-  \leq - \frac{1}{K_\delta} ,$$
  for $|\xi| \leq \delta$. And we obtain the  equivalence $ \lambda_+ + \lambda _-  \approx -1.$
  
   \item[{\rm (iv)}.] We observe that we have chosen $\eta > 0$ in \eqref{def-eta} such that 
   \begin{equation*}
 \frac{1}{2}\leq \sqrt{1 - 4 \log (1+ |\xi |^2 ) (1+\log (1+|\xi |^2))^2} \leq 1  
 \end{equation*}
for all $ |\xi | \leq \eta.  $
Thus, one can deduce
$$  \frac{1}{2 K_\delta }\leq \frac{1}{2(1+ \log (1+|\xi |^2))}\leq \frac{\sqrt{1 - 4 \log (1+ |\xi |^2 ) (1+\log (1+|\xi |^2))^2}}{1+ \log (1+|\xi |^2)} \leq \frac{1}{1+ \log (1+|\xi |^2)} \leq 1 .$$
This estimate shows the desired equivalence $ \lambda_+ - \lambda _-  \approx 1.$
\end{itemize}

\hfill
$\Box$

In the next subsection to use Lemma \ref{lemadeequivalencias} we work on the zone $\{|\xi| \leq \eta\}$, where $\eta$ is given in \eqref{def-eta}.


\subsection{Estimates on the low frequency zone  $ |\xi |  \leq \delta  $}

\textbf{(i) Estimates on the low frequency zone  $  |\xi |  \leq \eta  $}: 
\vspace{0.2cm}

We remember that $\eta$ is defined in \eqref{def-eta}. In this case, the characteristics roots $\lambda _{\pm }$ are real-valued, and the solution of \eqref{Probl-Fourier} is explicitly given by
\begin{align}\label{solutionlowfrequency}
\hat{u}(t, \xi) =  \frac{\lambda _- \hat{u}_0(\xi) - \hat{u}_1(\xi)}{\lambda _- - \lambda _+} e^{t \lambda _+} + \frac{\hat{u}_1(\xi) - \lambda _+ \hat{u}_0(\xi)  }{\lambda _- - \lambda _+} e ^{t \lambda _-}.
\end{align}
We observe that 
\begin{align*}
&\lambda _- = -\log (1+|\xi |^2 ) (1+ \log (1+ |\xi |^2)) - (1+\log (1+|\xi |^2 ))\lambda _-^2,  \\
&\lambda _+ = -\log (1+|\xi |^2 ) (1+ \log (1+ |\xi |^2)) - (1+\log (1+|\xi |^2 ))\lambda _+^2. 
\end{align*}
Therefore we can rewrite $ \hat{u}(t, \xi) $ as follows
\begin{align}
\hat{u}(t, \xi) =  e^{-t\log(1+|\xi|^2)(1+\log(1+|\xi |^2))}\left ( H_1(t, \xi) +H_2(t, \xi)  \right ),
\end{align}
where 
\begin{align*}
&H_1(t,\xi) = \frac{\lambda _- \hat{u}_0(\xi) - \hat{u}_1(\xi)}{\lambda _- - \lambda _+} e^{-t(1+\log(1+|\xi|^2))\lambda _+^2 }  ,\\
&H_2(t,\xi)= \frac{\hat{u}_1(\xi) - \lambda _+ \hat{u}_0(\xi)  }{\lambda _- - \lambda _+} e ^{-t(1+\log(1+|\xi|^2))\lambda _-^2}.
\end{align*}
We can also use the Chill-Haraux \cite{C-H-01} idea that has also been used in \cite{I-S2016} to decompose $H_1(t,\xi)$ as 
\begin{align*}
H_1(t,\xi) &= \hat{u}_0 + \hat{u}_1 + \frac{\lambda _+ - \lambda _-}{\lambda _- - \lambda _+} \hat{u}_0 -  \frac{\lambda _+ - \lambda _-}{\lambda _+ - \lambda _-} \hat{u}_1 + H_1(t,\xi)\\
&= \hat{u}_0 + \hat{u}_1 - \frac{ \lambda _+}{\lambda _+ - \lambda _-} \hat{u}_0 + \frac{\lambda _- \hat{u}_0 \left (  e^{-t(1+\log(1+|\xi |^2)) \lambda _+ ^2 } -1 \right ) }{\lambda _- - \lambda _+}  \\
&+  \frac{ \hat{u}_1 \left ( e^{-t(1+\log(1+|\xi |^2)) \lambda _+^2 } - (\lambda _+ - \lambda _-) \right ) }{\lambda _+ - \lambda _-}.
\end{align*}
By combining the last expression together with the decomposition $\hat{u}_j(\xi) = A_j(\xi) - i B_j(\xi) +P_j $ for initial data given in \eqref{decompo} we can get the following expression for $\hat{u}(t,\xi)$   which holds for $|\xi| \leq \eta$: 
\begin{align}\label{expansaoparau}
\hat{u}(t, \xi) &=  e^{-t\log(1+|\xi|^2)(1+\log(1+|\xi |^2))} \Big( A_0(\xi ) - i B_0(\xi) +P_0 + A_1(\xi ) -iB_1(\xi) +P_1 \Big) \nonumber \\
&- e^{-t\log(1+|\xi|^2)(1+\log(1+|\xi |^2))}  \frac{ \lambda _+}{\lambda _+ - \lambda _-} \hat{u}_0 + e^{-t\log(1+|\xi|^2)(1+\log(1+|\xi |^2))}  \frac{\lambda _- \hat{u}_0 \left (  e^{-t(1+\log(1+|\xi |^2)) \lambda _+^2 } -1  \right ) }{\lambda _- - \lambda _+} \nonumber \\ 
&+ e^{-t\log(1+|\xi|^2)\big(1+\log(1+|\xi |^2)\big)} \frac{ \hat{u}_1 \left ( e^{-t(1+\log(1+|\xi |^2) ) \lambda _+^2 } - (\lambda _+ - \lambda _-) \right ) }{\lambda _+ - \lambda _-} \nonumber \\
& + e^{-t\log(1+|\xi|^2)\big(1+\log(1+|\xi |^2)\big)} H_2(t, \xi) .
\end{align}
Our main goal in this subsection is to introduce an asymptotic profile of the solution $\hat{u}(t,\xi)$ in the low frequency region as $t \rightarrow \infty$ in a simple form as 
\begin{align}\label{defvarphi_1}
\varphi _1 (t,\xi) := (P_0 + P_1) e^{-t\log(1+|\xi|^2)(1+\log(1+|\xi |^2))}  .
\end{align}
For this purpose,  it is necessary to find suitable estimates for the other six terms of the expression \eqref{expansaoparau} defined by the functions
\begin{align*}
&F_1(t,\xi) = e^{-t\log(1+|\xi|^2)(1+\log(1+|\xi |^2))} \left ( A_0(\xi ) - i B_0(\xi)    \right ), \\
&F_2(t,\xi) = e^{-t\log(1+|\xi|^2)(1+\log(1+|\xi |^2))} \left ( A_1(\xi ) - i B_1 (\xi)    \right ), \\
&F_3(t, \xi) = -  e^{-t\log(1+|\xi|^2)(1+\log(1+|\xi |^2))}  \frac{ \lambda _+}{\lambda _+ - \lambda _-} \hat{u}_0, \\
&F_4 (t,\xi) = e^{-t\log(1+|\xi|^2)(1+\log(1+|\xi |^2))}  \frac{\lambda _- \hat{u}_0 \left (  e^{-t(1+\log(1+|\xi |^2)) \lambda _+^2} -1  \right ) }{\lambda _- - \lambda _+} ,\\
&F_5 (t,\xi) = e^{-t\log(1+|\xi|^2)(1+\log(1+|\xi |^2))} \frac{ \hat{u}_1 \left ( e^{-t(1+\log(1+|\xi |^2) ) \lambda _+^2 } - (\lambda _+ - \lambda _-) \right ) }{\lambda _+ - \lambda _-} , \\
&F_6 (t,\xi) = e^{-t\log(1+|\xi|^2)(1+\log(1+|\xi |^2))} H_2(t, \xi) .
\end{align*}
Therefore, from \eqref{expansaoparau} and \eqref{defvarphi_1}, for $|\xi|\leq \eta $ we have
\begin{align}\label{expansaobaixafreq}
\hat{u}(t,\xi) - \varphi_1(t,\xi) = \sum _{j=1}^{6} F_j(t,\xi).
\end{align}

In order to estimate the difference given by \eqref{expansaobaixafreq} on the zone of low frequency $\{|\xi|\leq \eta\}$ we shall develop the next computations based on Lemmas \ref{lema2.6} and \ref{lemadeequivalencias}.

Now, we first observe that 
$$ 1 \leq 1+\log (1+ |\xi |^2 ) \leq  1+\log (1+ \eta^2 )  =: k_\eta , \quad |\xi| \leq \eta.$$
Then, for $j=0,1 $ by using Lemma \ref{lema2.6} with $\kappa=1$ one has 
\begin{align*}
\int _{|\xi |\leq \eta } e^{-2t\log(1+|\xi|^2)(1+\log(1+|\xi |^2))} |A_j (\xi) - i B_j(\xi) |^2 d \xi & \leq  \int _{|\xi |\leq \eta } e^{-2t\log(1+|\xi|^2) } |A_j (\xi) - i B_j(\xi) |^2 d \xi \\
 & = \int _{|\xi |\leq \eta } (1+ |\xi|^2)^{-2t}  |A_j (\xi) - i B_j(\xi) |^2 d \xi \\
&\leq (L+M)^2 \| u_j \|_{1,1} ^2 \int _{|\xi |\leq \eta } (1+ |\xi|^2)^{-2t} | \xi |^2  d \xi \\
&= \omega _n (L+M)^2 \| u_j \|_{1,1} ^2 \int _{0 }^{\eta} (1+ r^2)^{-2t} r^{n+1} d \xi\\
&\leq  \omega _n (L+M)^2 \| u_j \|_{1,1} ^2 \int _{0 }^{1} (1+ r^2)^{-2t} r^{n+1} d \xi \\
&\leq C \| u_j \|_{1,1} ^2 t^{-\frac{n+2}{2}}
\end{align*}
for $t \gg 1$ due to Lemma \ref{general-p}. Consequently, for $t \gg 1$ we have 
$$ \int _{|\xi |\leq \eta } |F_1 (t, \xi )|^2 d\xi  \leq C \| u_0 \|_{1,1} ^2 t^{-\frac{n+2}{2}} \quad \mbox{ and } \quad \int _{|\xi |\leq \eta } |F_2 (t, \xi ) |^2  d\xi  \leq C \| u_1 \|_{1,1} ^2 t^{-\frac{n+2}{2}}.$$
 
In order to get an estimate on the function $F_3(t,\xi)$ we also use the simple inequality $\log (1+r^2) \leq r^2  $ for all $r\geq 0$. Then, by using the items {\rm (i)} and {\rm (iv)} of Lemma \ref{lema2.6} we obtain that 
\begin{align*}
\int _{|\xi |\leq \eta } | F_3(t,\xi ) |^2 d \xi &= \int _{|\xi |\leq \eta } e^{-2t\log(1+|\xi|^2)(1+\log(1+|\xi |^2))} \left (  \frac{ \lambda _+}{\lambda _+ - \lambda _-} \right )^2 |\hat{u}_0 |^2 d \xi  \\
&\leq C \int _{|\xi |\leq \eta } e^{-2t\log(1+|\xi|^2)}   \log ^2(1+|\xi |^2 ) |\hat{u}_0 |^2 d \xi \\
&\leq C \omega _n \| u_0 \|_1^2 \int _{0}^{\eta } (1+r^2)^{-2t} r^{n+3} dr \\
&\leq C \omega _n \| u_0 \|_1^2 \int _{0}^{1 } (1+r^2)^{-t} r^{n+3} dr \\
&\leq C \omega _n \| u_0 \|_1^2 t^{ - \frac{n+4}{2}}, \quad t \gg 1,
\end{align*}
where one has just used that $|\hat{u}_0(\xi)| \leq  \| u_0 \|_{1}$ for all $ \xi \in \R^n$ and the simple fact that 
\[\left\vert\frac{ \lambda_{\pm}}{\lambda _+ - \lambda _-}\right\vert \leq C, \quad  \vert\xi\vert \leq \eta, \]
with $C > 0$ a constant.
\vspace{0.2cm}

To estimate the function $F_4(t,\xi)$ we need the elementary inequality  \begin{equation} \label{e-a}
|e ^{-a} - 1 | \leq |a|, \quad   a \geq 0. 
\end{equation}
Thus, we have 
\begin{align*}
\int _{|\xi |\leq \eta } | F_4(t,\xi ) |^2 d \xi & = \int _{|\xi | \leq \eta } e^{-2t\log(1+|\xi|^2)(1+\log(1+|\xi |^2))}  \left ( \frac{\lambda _-   }{\lambda _- - \lambda _+} \right )^2 \left (  e^{-t(1+\log(1+|\xi |^2)) \lambda _+^2 } -1 \right ) ^2| \hat{u}_0 |^2 d \xi \\ 
&\leq C \int _{|\xi | \leq \eta } e^{-2t\log(1+|\xi|^2)(1+\log(1+|\xi |^2))}   \left (  e^{-t(1+\log(1+|\xi |^2)) \lambda _+^2 } -1 \right ) ^2| \hat{u}_0 |^2 d \xi \\
&\leq C t^{2} \int _{|\xi | \leq \eta } e^{-2t\log(1+|\xi|^2)(1+\log(1+|\xi |^2))} (1+\log (1+|\xi |^2))^2 \lambda _+^4   | \hat{u}_0 |^2 d \xi \\
&\leq C  t^{2} \| u_0 \|_{1}^2 \int _{|\xi | \leq \eta } (1+|\xi |^2 )^{-2t}  \log ^4 (1+|\xi |^2)    d \xi \\
&\leq  C  t^{2} \| u_0 \|_{1}^2  \int _{|\xi | \leq \eta } (1+|\xi |^2 )^{-2t}  |\xi |^8   d \xi \\
&\leq  \omega _n C t^{2} \| u_0 \|_{1}^2 t^{-\frac{n+8}{2}} \\
&= C_n  \| u_0 \|_{1}^2 t^{-\frac{n+4}{2}}, \quad t \gg 1
\end{align*}
because  of $1 + \log(1+|\xi|^2) \leq 1+\log 2$, for $ |\xi|\leq \eta <1$, where we also used the fact that $\vert \lambda_{+}\vert \leq C\log(1+\vert\xi\vert^{2})$ for $\vert\xi\vert \leq \eta$. The constant $C_n>0$ depends only on $n$. 
\vspace{0.2cm}

Now, let $ D = 1- 4 \log (1+|\xi |^2) (1+ \log (1+|\xi |^2))^2 $. Then we observe that $D\geq 0$ for $|\xi| \leq \eta$ and 
\begin{align*}
& 1- (\lambda_+ - \lambda_- )  = \frac{1}{1+ \log (1+|\xi |^2)}  \frac{ 2 \log (1+ |\xi |^2) + \log ^2(1+ |\xi |^2)  + 4 \log (1+|\xi |^2) (1+ \log (1+|\xi |^2))^2}{1+ \log (1+|\xi |^2) + \sqrt{D}} . 
\end{align*}
In particular,  $ 1- (\lambda_+ - \lambda_- ) $ is positive and there exists a constant $ C _\eta >0  $ such that  
$$ | 1- (\lambda_+ - \lambda_- ) | \leq  C_\eta \log (1+ |\xi |^2) $$
for all $|\xi| \leq \eta$, where $\eta$ is defined in \eqref{def-eta}.

At this point we use the inequalities \eqref{e-a} and  $(a-b)^2 \leq 2a^2 +2b^2$ to get
the next estimate to the remainder function $F_5(t,\xi)$.
 \begin{align*}
\int _{|\xi |\leq \eta }& | F_5(t,\xi ) |^2 d \xi  = \int _{|\xi | \leq \eta} e^{-2t\log(1+|\xi|^2)(1+\log(1+|\xi |^2))} \frac{  \left ( e^{-t(1+\log(1+|\xi |^2) ) \lambda _+^2 } - (\lambda _+ - \lambda _-) \right )^2 }{(\lambda _+ - \lambda _-)^2} |\hat{u}_1 | d \xi   \\
& \leq  2 \int _{|\xi | \leq \eta} e^{-2t\log(1+|\xi|^2)(1+\log(1+|\xi |^2))} \frac{  \left ( e^{-t(1+\log(1+|\xi |^2) ) \lambda _+^2 } - 1 \right )^2 }{(\lambda _+ - \lambda _-)^2} |\hat{u}_1 |^2 d \xi \\
&+ 2 \int _{|\xi | \leq \eta} e^{-2t\log(1+|\xi|^2)(1+\log(1+|\xi |^2))} \frac{  \left ( 1 - (\lambda _+ - \lambda _-) \right )^2 }{(\lambda _+ - \lambda _-)^2} |\hat{u}_1 |^2 d \xi \\
&\leq C  \int _{|\xi | \leq \eta} e^{-2t\log(1+|\xi|^2)(1+\log(1+|\xi |^2))}   \frac{  t^2(1+\log(1+|\xi |^2) )^2 \lambda _+^4 }{(\lambda _+ - \lambda _-)^2} |\hat{u}_1 | ^2d \xi \\ 
&+ 2 C_\eta \int _{|\xi | \leq \eta} e^{-2t\log(1+|\xi|^2)(1+\log(1+|\xi |^2))} \log ^2 (1+ |\xi |^2) |\hat{u}_1 |^2 d \xi  \\
&\leq C t^2  \int _{|\xi | \leq \eta} e^{-t\log(1+|\xi|^2) }  \log ^4(1+|\xi |^2) |\hat{u}_1 |^2 d \xi  + 2 C_\eta \int _{|\xi | \leq \eta} e^{-t\log(1+|\xi|^2) } \log ^2 (1+ |\xi |^2) |\hat{u}_1 |^2 d \xi \\
&\leq  C t^2 \| u_1 \|_1^2  \int _{|\xi | \leq 1} (1+|\xi |^2)^{-t }   |\xi |^8  d \xi +  2 C_\eta \|u_1 \|_1^2 \int _{|\xi | \leq 1} (1+|\xi|^2)^{-t} | \xi |^4  d \xi \\
&\leq C_n  \| u_1 \|_1^2 t^{-\frac{n+4}{2}} + 2 C_{\eta, n } \|u_1 \|_1^2  t^{-\frac{n+4}{2}}, \quad t \gg 1.
\end{align*}

Finally, by {\rm (ii)} of Lemma \ref{lemadeequivalencias} one has $ \lambda _- \approx -1  $ on the region $|\xi| \leq \eta$, so that there exists constants $c_1, c_2 >0 $ such that  $ c_1 \leq 2(1+\log(1+|\xi|^2))\lambda _-^2 \leq c_2$ whenever $|\xi | \leq \eta$. Therefore, it follows that
\begin{align*}
\int _{|\xi |\leq \eta } | F_6(t,\xi ) |^2 d \xi &= \int _{|\xi |\leq \eta }e^{-2t\log(1+|\xi|^2)(1+\log(1+|\xi |^2))} H_2^2(t, \xi) d \xi \\
&\leq  \int _{|\xi |\leq \eta }e^{-2t\log(1+|\xi|^2)(1+\log(1+|\xi |^2))}  e ^{-2t(1+\log(1+|\xi|^2))\lambda _-^2} \frac{1 }{(\lambda _- - \lambda _+)^2} |\hat{ u}_1 | d \xi  \\
&+ \int _{|\xi |\leq \eta }e^{-2t\log(1+|\xi|^2)(1+\log(1+|\xi |^2))} e ^{-2t(1+\log(1+|\xi|^2))\lambda _-^2} \left ( \frac{ \lambda _+  }{\lambda _- - \lambda _+} \right )^2  |\hat{u}_0 | d \xi \\
& \leq C e ^{-c_2t} \int _{|\xi |\leq \eta }e^{-t\log(1+|\xi|^2) }  | \hat{u}_1 |^2 d \xi    +  C e ^{-c_2t} \int _{|\xi |\leq \eta }e^{-t\log(1+|\xi|^2)}  \log ^2 (1+ |\xi |^2 ) | \hat{u}_0 |^2 d \xi   \\
& \leq C e ^{-c_2t} \| u_1 \|_1^2 \int _{|\xi |\leq \eta } (1+|\xi |^2)^{-t}  d \xi    +  C e ^{-c_2t} \| u_0 \|_1^2 \int _{|\xi |\leq \eta }e^{-t\log(1+|\xi|^2)} |\xi |^4  d \xi   \\
& \leq C \| u_1 \|_1^2 t^{-\frac{n}{2}}  e ^{-c_2t}  +  C  \| u_0 \|_1^2 t^{-\frac{n+4}{2}}  e ^{-c_2t} , \qquad t \gg 1.
\end{align*}

By combining the above estimates for $F_j(t,\xi), \; j=1, \cdots , 6,$ with equation \eqref{expansaobaixafreq}, we obtain that the solution $\hat{u}(t,\xi)$ given by \eqref{expansaoparau} has the following asymptotic property.

\begin{pro}\label{proplowfrequency}
Let $ n\geq 1  $ and $ \eta >0 $ given by \eqref{def-eta}. For  $ (u_0, u_1) \in L^{1,1}({\bf R}^{n}) \times L^{1,1}({\bf R}^{n}) $   the solution  $ \hat{u} (t,\xi ) $ to problem \eqref{Probl-Fourier} satisfies 
\[\int _{|\xi | \leq \eta } | \hat{u}(t,\xi ) - \varphi_1(t,\xi) |^2 d \xi 
\leq C  \Big( \| u_0 \|_{1,1} ^2 t^{-\frac{n+2}{2}} + \| u_1 \|_{1,1} ^2 t^{-\frac{n+2}{2}} + \| u_0 \|_1^2 t^{ - \frac{n+4}{2}}\] 
\[+ \| u_1 \|_1^2 t^{-\frac{n+4}{2}} +  \|u_1 \|_1^2  t^{-\frac{n+2}{2}} + \| u_1 \|_1^2 t^{-\frac{n}{2}}  e ^{-c_2t}  +   \| u_0 \|_1^2 t^{-\frac{n+4}{2}}  e ^{-c_2t} \Big), \quad t\gg 1\]
where $\varphi_1(t,\xi) $ is defined in \eqref{defvarphi_1} and $C$ is a positive constant that depends only on $ \eta $ and $ n $.  
\end{pro} 
\begin{flushright}
$\square$ 
\end{flushright}


\textbf{(ii) Estimates on the low-middle frequency zone  $ \eta \leq |\xi |  \leq \delta  $}:
\vspace{0.2cm}

When one obtains sharp estimates on the region middle frequency zone $ \eta \leq |\xi | \leq \delta $ it should be noted that according to \eqref{condicao} the characteristics roots $\lambda_{\pm}$ are still real-valued. Then we can  rewrite the solution $ \hat{u}(t,\xi)$ as follows.
 \begin{align*}
 \hat{u}(t, \xi) &= e^{-\frac{t}{2 (1+\log(1+|\xi |^2))}} \left [ \cosh (c(\xi)t) \hat{u}_0 + \frac{\sinh(c(\xi)t)}{2(1+\log(1+|\xi |^2))c(\xi)} \hat{u}_0 + \frac{\sinh( c(\xi)t )}{c(\xi )} \hat{u}_1 \right ] 
 \end{align*}
where
$$ c(\xi):= \frac{ \sqrt{1-4 \log (1+|\xi |^2 ) (1+\log(1+|\xi |^2))^2}}{2 (1+\log (1+|\xi|^2))}.$$

We remember that $\delta$ is given in \eqref{condicao} and our choice for $ \eta $ is such that 
$$  \sqrt{1 - 4 \log (1+ |\xi |^2 ) (1+\log (1+|\xi |^2))^2}   \geq \frac{1}{2} $$
when $|\xi| \leq \eta$ (see \eqref{def-eta}), and this is the decreasing function on $\vert\xi\vert$. Thus, in the case for $\eta \leq | \xi | \leq  \delta$, one has 
 $$ c(\xi ) \leq \frac{1}{4 (1+\log (1+|\xi | ^2)) } .$$
Then, due to the fact that  $  \cosh a \leq e^a $ for all $ a \geq 0 $ we may estimate for $t >0$
\begin{align*}
\int _{\eta \leq |\xi | \leq \delta } e^{-\frac{t}{ 1+\log(1+|\xi |^2)}} \cosh ^2 (c(\xi)t) |\hat{u}_0 (\xi ) |^2 d\xi &\leq \int _{\eta  \leq |\xi | \leq \delta } e^{-\frac{t}{ 1+\log(1+|\xi |^2)}} e^{ 2 c(\xi)t} |\hat{u}_0 (\xi ) |^2 d\xi \\
& \leq \int _{\eta \leq |\xi | \leq \delta } e^{-\frac{t}{ 1+\log(1+|\xi |^2)}} e^{ \frac{t}{2 (1+\log(1+|\xi|^2))}} |\hat{u}_0 (\xi ) |^2 d\xi \\
&= \int _{\eta \leq |\xi | \leq \delta }  e^{ - \frac{t}{2 (1+\log(1+|\xi|^2))}} |\hat{u}_0 (\xi ) |^2 d\xi \\
&\leq e ^{-C_\delta t}  \int _{\eta \leq |\xi | \leq \delta }  |\hat{u}_0 (\xi ) |^2 d\xi \\
& \leq e ^{-C_\delta t} \| u_0 \|_2^2 ,
\end{align*}
with $C_{\delta}=\frac{1}{2(1+\log(1+\delta^2))}$.
\vspace{0.2cm}

By combining Lemma \ref{lemmahiperbolicsine} with ideas  used to get the previous estimate, we can  obtain $L^2$-estimates for the other two terms of the solution $\hat{u}(t,\xi)$ in the middle frequency zone. In fact, we have 
\begin{align*}
\int _{\eta \leq |\xi | \leq \delta } e^{-\frac{t}{ 1+\log(1+|\xi |^2)}} \frac{\sinh ^2( c(\xi)t )}{c(\xi )^2} | \hat{u}_1 |^2 d \xi &\leq \int _{\eta \leq |\xi | \leq \delta } e^{-\frac{t}{ 1+\log(1+|\xi |^2)}} t^2 e^{2c(\xi)t} | \hat{u}_1 |^2 d \xi \\
 &\leq t^{2}e^{-C_\delta t } \| u_1 \| _2^2, 
\end{align*}
and 
\begin{align*}
 \int _{\eta \leq |\xi | \leq \delta } e^{-\frac{t}{ 1+\log(1+|\xi |^2)}} \frac{\sinh ^2(c(\xi)t)}{4(1+\log(1+|\xi |^2))^2c(\xi)^2 } | \hat{u}_0|^2d\xi &\leq t^2  \int _{\eta \leq |\xi | \leq \delta }  \frac{ e^{-\frac{t}{ 1+\log(1+|\xi |^2)}} }{4(1+\log(1+|\xi |^2))^2 } e^{2c(\xi)t } | \hat{u}_0|^2 d\xi \\
& \leq C_\eta t^2 e^{-C_\delta t} \| u _0 \|_2^2 
\end{align*}
for all $t>0$. Therefore, the following estimate on the middle frequency zone holds for $t>0$
\begin{align}
 \int _{\eta \leq |\xi | \leq \delta } | \hat{u} (t,\xi) | ^2 d \xi & \leq \int _{\eta \leq |\xi | \leq \delta } e^{-\frac{t}{ 1+\log(1+|\xi |^2)}} \cosh ^2 (c(\xi)t) |\hat{u}_0 (\xi ) |^2 d\xi \nonumber \\ 
 &+  \int _{\eta \leq |\xi | \leq \delta } e^{-\frac{t}{ 1+\log(1+|\xi |^2)}} \frac{\sinh ^2( c(\xi)t )}{c(\xi )^2} | \hat{u}_1 |^2 d \xi \nonumber \\ 
 &+  \int _{\eta \leq |\xi | \leq \delta } e^{-\frac{t}{ 1+\log(1+|\xi |^2)}} \frac{\sinh ^2(c(\xi)t)}{4(1+\log(1+|\xi |^2))^2c(\xi)^2 } | \hat{u}_0|^2d\xi \nonumber \\
 & \leq e ^{-C_\delta t} \| u_0 \|_2^2 + C_\eta t^2 e^{-C_\delta t} \| u _0 \|_2^2 +  t^{2}e^{-C_\delta t } \| u_1 \| _2^2. 
 \end{align}
By summarizing the above estimates in the middle frequency zone one can conclude the next result.  
\begin{lem}\label{lemmadecayulow-middle}
Let $ n\geq 1  $ and $ u_0, u_1 \in L^2({\bf R}^{n})$. Then the solution  $ \hat{u} (t,\xi ) $ to the problem \eqref{Probl-Fourier} satisfies 
\begin{equation}
 \int _{\eta \leq |\xi | \leq \delta } | \hat{u} (t,\xi) | ^2 d \xi  \leq e ^{-C_\delta t} \| u_0 \|_2^2 + C_\eta t^2 e^{-C_\delta t} \| u _0 \|_2^2 +  t^{2}e^{-C_\delta t } \| u_1 \| _2^2
\end{equation}
for $ t>0 $.  \begin{flushright}
$\square$ 
\end{flushright}
\end{lem}

\vspace{-0.7cm}

\subsection{Estimates on the high frequency zone  $ |\xi | \geq \delta  $}

On the zone of high frequency the characteristics roots are complex-valued (see \ref{condicao}), and are given by
$$ \lambda _\pm =-a(\xi) \pm i b(\xi),  $$
where \begin{equation}\label{definition a and b}
a(\xi ) = \frac{1}{2(1+\log (1+|\xi|^2))} \text{\; and\; } b(\xi ) = \frac{\sqrt{4  \log(1+|\xi|^2)(1+\log(1+|\xi|^2))^2-1} }{2 (1+\log (1+|\xi|^2))} .
\end{equation}
Then the solution $ \hat{u }(t,\xi)$ to problem \eqref{Probl-Fourier} in the high frequency zone is explicitly given by
\begin{equation}\label{solutioninhighfrequency}
\hat{u}(t, \xi) =  e^{-a(\xi)t} \cos (b(\xi)t) \hat{u}_0 + \frac{a(\xi) }{b(\xi)} e^{-a(\xi)t}  \sin (b(\xi)t) \hat{u}_0 + \frac{1}{b(\xi)} e^{-a(\xi)t} \sin (b(\xi)t) \hat{u}_1 .
\end{equation}

\textbf{(i)  Estimate on the high-middle frequency zone $ \delta \leq |\xi | \leq \sqrt{e-1}  $}. 
\vspace{0.2cm} 

In this region, the function $ a(\xi) $ is equivalent to a constant, that is $ \frac{1}{4} \leq a(\xi ) \leq \frac{1}{2} .$

Moreover, we can see that $\displaystyle{\frac{1}{b(\xi)}}$ converges to  $+\infty $ when $|\xi | \rightarrow \delta ^+ $ according to \eqref{condicao}. However, we remember that $ \sin a \leq a  $ for all $ a \geq 0 $. Thus 
$$\sin (b(\xi)t) \leq b(\xi) t$$
for all $ \xi \in {\bf R}^{n}$ and $ t \geq 0 $.
By combining these properties together with the solution formula \eqref{solutioninhighfrequency} one can obtain the following estimate for $ t >0 $, which implies the exponential decay in such region.  
\begin{equation}
\int _{\delta \leq |\xi | \leq \sqrt{e-1}} |\hat{u}(t,\xi) |^2 d\xi \leq e^{-\frac{t}{2}} \| u_0 \|_2^2+ \frac{1}{4} t^2e^{-\frac{t}{2}} \| u_0 \|_2^2 + t^2e^{-\frac{t}{2}} \| u_1 \|_2^2 .
\end{equation}

\textbf{(ii)  Estimate on the high frequency zone $  |\xi | \geq \sqrt{e-1}  $} 
\vspace{0.2cm} 

The estimates on this zone are more delicate and the derivation is one of essential contributions in our paper. We first need to rewrite the solution formula given by \eqref{solutioninhighfrequency} into a more suitable expression. 

First we observe that for $ |\xi | \geq  \delta  $, in particular, for $|\xi | \geq \sqrt{e-1}$, it holds that 
$$ b(\xi) \leq \sqrt{\log (1+|\xi|^2)} .$$
Then the mean value theorem implies,  for $|\xi| \geq \sqrt{e-1}$, that  
$$ \cos( b(\xi)t )   = \cos (\sqrt{\log (1+|\xi|^2)}t) -\sin(\theta (t, \xi)) \left [ b(\xi) - \sqrt{\log(1+|\xi|^2)} \right ]t, $$
with $ \theta (t, \xi) = \alpha b(\xi) t  +(1-\alpha ) \sqrt{\log(1+|\xi|^2)} t $ for some $ \alpha \in (0,1)$. 
\newline
Similarly, 
$$ \sin( b(\xi)t )   = \sin (\sqrt{\log (1+|\xi|^2)}t) +\cos(\eta (t, \xi)) \left [ b(\xi) - \sqrt{\log(1+|\xi|^2)} \right ]t, $$
with $ \eta (t, \xi) = \gamma b(\xi) t  +(1-\gamma ) \sqrt{\log(1+|\xi|^2)} t $ for some $ \gamma \in (0,1)$.
\vspace{0.2cm}

Thus, one can rewrite $\hat{u}(t, \xi) $ as follows:
\begin{align}\label{u-high}
 \hat{u}(t, \xi) &= e^{-a(\xi)t} \cos (\sqrt{\log (1+|\xi|^2)} t ) \hat{u}_0 + t e^{-a(\xi)t} \sin (\theta(\xi,t) \left [  \sqrt{\log (1+|\xi|^2)}  - b(\xi)\right ] \hat{u}_0 \nonumber  \\ 
&+ \frac{a(\xi) }{b(\xi)} e^{-a(\xi)t}  \sin (b(\xi)t) \hat{u}_0 + 
\frac{1}{b(\xi)} e^{-a(\xi)t} \sin (\sqrt{\log (1+|\xi|^2)} t) \hat{u}_1  \nonumber \\ 
&+   te^{-a(\xi)t} \frac{1}{b(\xi)} \cos(\eta(\xi,t)) \left [ b(\xi) - \sqrt{\log(1+|\xi |^2)} \right ] \hat{u}_1. 
\end{align}
We introduce  an important function to be the asymptotic profile on the zone of high frequency for the solution $\hat{u}(t,\xi)$ given by \eqref{u-high} as follows
\begin{align}\label{defvarphi_2}
\varphi _2(t, \xi) :=e^{-\frac{t}{2 \log (1+|\xi |^2 )}} \frac{\sin (\sqrt{\log (1+|\xi |^2 )}t )}{\sqrt{\log (1+|\xi |^2)}} \hat{u}_1(\xi) +  e^{-\frac{t}{2 \log (1+|\xi |^2 )}} \cos (\sqrt{\log (1+|\xi |^2 )}t ) \hat{u}_0(\xi).
\end{align}
In the following part, one will prove that the function $\varphi _2(t, \xi)$ is  asymptotic profile for the solution $\hat{u}(t,\xi)$ in the high frequency region $|\xi | \geq \sqrt{e-1}$. Then we obtain the following difference between the solution and the profile
\begin{align}\label{u-perfil}
\hat{u}(t, \xi)  - \varphi _2 (t,\xi) &= \left ( e^{-a(\xi)t} - e^{-\frac{t}{2 \log (1+|\xi |^2 )}}  \right )\cos (\sqrt{\log (1+|\xi|^2)} t ) \hat{u}_0 +  \frac{a(\xi) }{b(\xi)} e^{-a(\xi)t}  \sin (b(\xi)t) \hat{u}_0  \nonumber \\
&   + t e^{-a(\xi)t} \sin (\theta(\xi,t)) \left [  \sqrt{\log (1+|\xi|^2)}  - b(\xi)\right ] \hat{u}_0  \nonumber \\ 
&+   e^{-a(\xi)t} \left ( \frac{1}{b(\xi)}  -  \frac{1 }{\sqrt{\log (1+|\xi |^2)}}  \right )  \sin (\sqrt{\log (1+|\xi|^2)} t) \hat{u}_1  \nonumber  \\ 
&+ \frac{1 }{\sqrt{\log (1+|\xi |^2)}} \left (  e^{-a(\xi)t} -   e^{-\frac{t}{2 \log (1+|\xi |^2 )}} \right )  \sin (\sqrt{\log (1+|\xi|^2)} t) \hat{u}_1  \nonumber  \\ 
&+   te^{-a(\xi)t} \frac{1}{b(\xi)} \cos(\eta(\xi,t)) \left [ b(\xi) - \sqrt{\log(1+|\xi |^2)} \right ] \hat{u}_1 .
\end{align}

Then, the following functions 
\begin{align*}
&G_1(t, \xi) = \left ( e^{-a(\xi)t} - e^{-\frac{t}{2 \log (1+|\xi |^2 )}}  \right )\cos (\sqrt{\log (1+|\xi|^2)} t ) \hat{u}_0 , \\ 
&G_2(t, \xi) = t e^{-a(\xi)t} \sin (\theta(\xi,t)) \left [  \sqrt{\log (1+|\xi|^2)}  - b(\xi)\right ] \hat{u}_0  ,\\ 
&G_3(t, \xi) =  \frac{a(\xi) }{b(\xi)} e^{-a(\xi)t}  \sin (b(\xi)t) \hat{u}_0 , \\
 &  G_4(t, \xi) =  e^{-a(\xi)t} \left ( \frac{1}{b(\xi)}  -  \frac{1 }{\sqrt{\log (1+|\xi |^2)}}  \right )  \sin (\sqrt{\log (1+|\xi|^2)} t) \hat{u}_1  ,\\
  &  G_5(t, \xi) =  \frac{1 }{\sqrt{\log (1+|\xi |^2)}} \left (  e^{-a(\xi)t} -   e^{-\frac{t}{2 \log (1+|\xi |^2 )}} \right )  \sin (\sqrt{\log (1+|\xi|^2)} t) \hat{u}_1 , \\
&G_6(t, \xi) = te^{-a(\xi)t} \frac{1}{b(\xi)} \cos(\eta(\xi,t)) \left [ b(\xi) - \sqrt{\log(1+|\xi |^2)} \right ] \hat{u}_1
\end{align*}
are the remainder terms that appear in \eqref{u-perfil}.
\vspace{0.2cm}

From now, let us estimates these $6$-remainders in the following lines. 

We note that on the region such that $| \xi | \geq \sqrt{e-1}$ one has $1+ \log (1+|\xi |^2 ) \leq 2 \log (1+|\xi |^2 ) $. Also, by Lemma \ref{LemmaTecnico1} with $c=1$ and $a = -1$ one has 
\begin{equation}\label{decay-nu}
 \frac{e^{\frac{-t}{1+\log(1+|\xi |^2 ) }}}{(1+\log(1+|\xi |^2 ))^\nu } \leq C t^{-\nu} , \quad t>0, \;\; \xi \in {\bf R}^{n}, 
 \end{equation}
for a fixed $\nu >0$. The above two inequalities will be used to get the next series of estimates for the functions $G_j(t,\xi),\, j=1 \cdots, 6$.

 The first estimate is concerned with the function $G_1(t,\xi)$. 
\begin{align*}
\int _{|\xi |\geq \sqrt{e-1} } |G_1(t, \xi)|^2 d \xi &= \int _{|\xi |\geq \sqrt{e-1}} \left ( e^{-a(\xi)t} - e^{-\frac{t}{2 \log (1+|\xi |^2 )}}  \right )^2 \cos ^2(\sqrt{\log (1+|\xi|^2)} t ) |\hat{u}_0| ^2  d \xi \\
&= \int _{|\xi |\geq \sqrt{e-1} } \hspace{-0.2cm} e^{\frac{-t}{1+\log (1+|\xi |^2 ) }} \left ( 1 - e^{\frac{-t}{2\log (1+|\xi |^2 ) (1+ \log (1+|\xi |^2 ))}} \right ) ^2 \cos ^2(\sqrt{\log (1+|\xi|^2)} t ) |\hat{u}_0| ^2  d \xi \\
&\leq \int _{|\xi |\geq \sqrt{e-1}} e^{\frac{-t}{1+\log (1+|\xi |^2 ) }} \left ( 1 - e^{\frac{-t}{2\log (1+|\xi |^2 ) (1+ \log (1+|\xi |^2 ))}} \right ) ^2 |\hat{u}_0| ^2  d \xi.
\end{align*}
The next inequalities conclude the estimate to the function $G_1(t,\xi)$.
\begin{align*}
\int _{|\xi |\geq \sqrt{e-1} } |G_1(t, \xi)|^2 d \xi &\leq t^2 \int _{|\xi |\geq \sqrt{e-1} }   \frac{e^{\frac{-t}{1+\log (1+|\xi |^2 ) }} }{4\log ^2 (1+|\xi |^2 ) (1+ \log (1+|\xi |^2 ))^2 }  |\hat{u}_0| ^2  d \xi \\
&\leq t^2 \int _{|\xi |\geq \sqrt{e-1}}   \frac{e^{\frac{-t}{1+\log (1+|\xi |^2 ) }}  }{ (1+ \log (1+|\xi |^2))^4 } |\hat{u}_0| ^2  d \xi \\
&= t^2 \int _{|\xi |\geq \sqrt{e-1}}   \frac{e^{\frac{-t}{1+\log (1+|\xi |^2 ) }}  }{ (1+ \log (1+|\xi |^2))^{5+l}} (1+ \log (1+|\xi |^2) )^{l+1} |\hat{u}_0| ^2  d \xi \\
&\leq C t^2 t^{-(l + 5 )} \int _{|\xi |\geq \sqrt{e-1} }    (1+ \log (1+|\xi |^2) )^{l+1} |\hat{u}_0| ^2  d \xi \\
& \leq C t^{-(l + 3 )} \| u_0 \| _{Y^{l+1} }^2, 
\end{align*}
for all $t>0$ and $l \geq 0$, where we have used the inequalities  \eqref{e-a}, \eqref{decay-nu} and the fact that $\log(1+\vert\xi\vert^{2}) \geq 1$ on the high frequency zone.
\vspace{0.2cm}

For $| \xi | \geq \delta $, we introduce the auxiliary function $R(t,\xi)$ defined by 
\begin{align}\label{defR(t,xi)}
R(t,\xi)= \sqrt{1 -\frac{1}{4 \log (1+|\xi |^2) (1+\log (1+|\xi |^2))^2}},
\end{align}
\vspace{0.2cm}
which is well defined  due to \eqref{condicao}. Additionally, one notes that for $|\xi | \geq \sqrt{e-1}$ we have the following estimate 
\begin{align*}
|\sqrt{\log(1+|\xi |^2)} - b(\xi) |& = \frac{1}{4 (1+\log (1+|\xi |^2))^2 \sqrt{\log(1+|\xi |^2)} \left ( 1+  R(t,\xi)  \right )} \\
& \leq \frac{1}{4 (1+\log (1+|\xi |^2))^2 \sqrt{\log(1+|\xi |^2)} }.  
\end{align*} 

Thus for $t >0 $ and $l \geq 0$ we get 
\begin{align*}
\int _{|\xi |\geq \sqrt{e-1} } |G_2(t, \xi)|^2 d \xi &= t^2 \int _{|\xi |\geq \sqrt{e-1} }   e^{-2a(\xi)t} \sin ^2  (\theta(\xi,t)) \left [  \sqrt{\log (1+|\xi|^2)}  - b(\xi)\right ]^2 | \hat{u}_0 |^2  d \xi \\
& \leq t^2 \int _{|\xi |\geq \sqrt{e-1} }       \frac{e^{-\frac{t}{1+ \log (1+|\xi |^2 )}} }{16 (1+\log (1+|\xi |^2))^4 \log(1+|\xi |^2) } | \hat{u}_0 |^2  d \xi \\ 
&\leq t^2 \int _{|\xi |\geq \sqrt{e-1} }       \frac{e^{-\frac{t}{1+ \log (1+|\xi |^2 )}} }{8 (1+\log (1+|\xi |^2))^5 } | \hat{u}_0 |^2  d \xi \\ 
&\leq t^2 \int _{|\xi |\geq \sqrt{e-1} }       \frac{e^{-\frac{t}{1+ \log (1+|\xi |^2 )}} }{8 (1+\log (1+|\xi |^2))^{6+l} }  (1+ \log (1+|\xi |^2 ))^{l+1} | \hat{u}_0 |^2  d \xi \\ 
&= C t^2 t^{-(6+l)} \int _{|\xi |\geq \sqrt{e-1} }   (1+ \log (1+|\xi |^2 ))^{l+1} | \hat{u}_0 |^2  d \xi \\ 
&\leq C t^{-(l+4)} \| u_0 \|_{Y^{l+1}}^2,
\end{align*}
where one has just used the fact that $ 1 + \log (1+|\xi|^2) \geq 2$ for $|\xi|\geq \sqrt{e-1} $ and \eqref{decay-nu}.

Another important property is that the function $\vert\xi\vert \mapsto \displaystyle{\frac{1}{4 \log (1+|\xi |^2) (1+\log (1+|\xi |^2))^2}}$ is decreasing, and it converges to zero as $ |\xi | \rightarrow +\infty $. Hence, it follows that
\begin{align}\label{ineq16/15}
1 \leq \frac{1}{ 1 -\frac{1}{4 \log (1+|\xi |^2) (1+\log (1+|\xi |^2))^2} } = \frac{1}{R(t,\xi)^2 } \leq \frac{16}{15} 
\end{align}
for  $ | \xi | \geq \sqrt{e-1} $.

From the above inequality one can obtain estimates of the $L^2-$norms of each functions $G_3(t, \cdot ) $, $G_4(t, \cdot ) $ and $ G_5(t, \cdot ) $ for $t >0$.  In fact,  \eqref{ineq16/15} implies that 
\begin{align*}
\int _{|\xi | \geq \sqrt{e-1}} |G_3(t,\xi) |^2 d \xi &= \int _{|\xi | \geq \sqrt{e-1}} \left ( \frac{a(\xi) }{b(\xi)} \right )^2 e^{-2a(\xi)t}  \sin ^2 (b(\xi)t) | \hat{u}_0|^2 d \xi \\
&=  \int _{|\xi | \geq \sqrt{e-1}}  \frac{ e^{-\frac{t}{1+ \log (1+|\xi |^2) } } \sin ^2 (b(\xi)t) | \hat{u}_0|^2 }{4 \log (1+|\xi |^2)(1+ \log (1+|\xi |^2))^2  R(t,\xi)^2  }     d \xi \\ 
&\leq \frac{16}{15}  \int _{|\xi | \geq \sqrt{e-1}}  \frac{ e^{-\frac{t}{1+ \log (1+|\xi |^2) } } }{4 \log (1+|\xi |^2)(1+ \log (1+|\xi |^2))^2   }     | \hat{u}_0|^2 d \xi \\ 
&\leq \frac{8}{15}  \int _{|\xi | \geq \sqrt{e-1}}  \frac{ e^{-\frac{t}{1+ \log (1+|\xi |^2) } } }{  (1+ \log (1+|\xi |^2))^3   }     | \hat{u}_0|^2 d \xi \\
&=\frac{8}{15}  \int _{|\xi | \geq \sqrt{e-1}}  \frac{ e^{-\frac{t}{1+ \log (1+|\xi |^2) } } }{  (1+ \log (1+|\xi |^2))^{4+l}   }  (1+ \log (1+|\xi |^2))^{l+1}   | \hat{u}_0|^2 d \xi \\
&\leq C t^{-(l+4)} \int _{|\xi | \geq \sqrt{e-1}}    (1+ \log (1+|\xi |^2))^{l+1}   | \hat{u}_0|^2 d \xi \\
&\leq C t^{-(4+l)} \| u_0 \|_{Y^{l+1}}^2,
\end{align*}
where we have used the fact that $ 4 \log (1+|\xi|^2) \geq 2+2 \log (1+|\xi|^2) $ for $|\xi | \geq \sqrt{e-1} $.
\vspace{0.2cm}


To get an estimate for the $L^2$-norm of $ G_4(t, \cdot ) $ we first observe the following identity: 
\begin{align*}
\frac{1}{b(\xi)} - \frac{1}{\sqrt{\log (1+|\xi |^2)}}
 = \frac{1}{ 4\log ^{3/2}(1+|\xi |^2) (1+ \log (1+|\xi |^2))^2  R(t,\xi) \left ( 1 + R(t,\xi)     \right )}  . 
\end{align*}
holds for $ |\xi | \geq \delta $, where $R(t,\xi)  $ is given by \eqref{defR(t,xi)}. 
\vspace{0.2cm}

By using the above identity, the estimate \eqref{ineq16/15} and the fact that $ 1+ \log (1+|\xi |^2) \leq 2  \log (1+|\xi |^2)  $ for $\vert\xi\vert \leq \sqrt{e-1}$, we can arrive at the following $L^2$-estimate to the function $ G_4(t,\cdot ) $: 
\begin{align*}
\int _{|\xi | \geq \sqrt{e-1}} |G_4(t, \xi) |^2 d \xi  &=  \int _{|\xi | \geq \sqrt{e-1}} e^{-2a(\xi)t} \left ( \frac{1}{b(\xi)}  -  \frac{1 }{\sqrt{\log (1+|\xi |^2)}}  \right )^2  \sin ^2 (\sqrt{\log (1+|\xi|^2)} t) |\hat{u}_1 |^2 d \xi  \\
&\leq \frac{1}{15} \int _{|\xi | \geq \sqrt{e-1}} \frac{e^{ - \frac{t}{(1+\log (1+|\xi |^2))}}}{\log ^3(1+|\xi |^2) (1+ \log (1+|\xi |^2 ))^4}  \sin ^2 (\sqrt{\log (1+|\xi|^2)} t) |\hat{u}_1 |^2 d \xi  \\
&\leq \frac{8}{15} \int _{|\xi | \geq \sqrt{e-1}} \frac{e^{ - \frac{t}{(1+\log (1+|\xi |^2))}}}{ (1+ \log (1+|\xi |^2 ))^7}   |\hat{u}_1 |^2 d \xi  \\
&\leq C t^{-(l+7)} \|u_1 \| _{Y^l}^2.
\end{align*}
Similarly to the previous estimate for $G_1(t, \cdot) $ one obtains
\begin{align*}
\int _{|\xi \geq \sqrt{e-1}} |G_5(t,\xi) |^2 d \xi &= \int _{|\xi \geq \sqrt{e-1}}  \frac{1 }{\log (1+|\xi |^2)} \left (  e^{- \frac{t}{1+ \log (1+|\xi |^2 )}} -   e^{-\frac{t}{2 \log (1+|\xi |^2 )}} \right ) ^2 \sin ^2 (\sqrt{\log (1+|\xi|^2)} t) | \hat{u}_1 |^2 d \xi  \\
&= \int _{|\xi \geq \sqrt{e-1}}  \frac{e^{ - \frac{t}{1+\log(1+|\xi |^2)}} }{\log (1+|\xi |^2)}  \left (  1 -   e^{-\frac{t}{2 \log (1+|\xi |^2 ) (1+ \log (1+|\xi |^2))}} \right ) ^2 \sin ^2 (\sqrt{\log (1+|\xi|^2)} t) | \hat{u}_1 |^2 d \xi  \\
&\leq t^2 \int _{|\xi \geq \sqrt{e-1}}  \frac{ e^{ - \frac{t}{1+\log(1+|\xi |^2)}} }{ \log (1+|\xi |^2)} \frac{1}{4 \log ^2 (1+|\xi |^2 ) (1+ \log (1+|\xi |^2))^2 }  | \hat{u}_1 |^2 d \xi  \\ 
&\leq 2  t^2 \int _{|\xi \geq \sqrt{e-1}}  \frac{e^{ - \frac{t}{1+\log(1+|\xi |^2)}}}{(1+ \log (1+|\xi |^2))^5 }  | \hat{u}_1 |^2 d \xi \\
& \leq C t^{-(l+3)} \| u_1 \|_{Y^l}^2.
\end{align*}
Finally, we observe that 
\begin{align*}
\frac{b(\xi) -\sqrt{\log (1+|\xi |^2 ) }}{b(\xi)} =&  \frac{ -1 }{4 \log (1+|\xi |^2 ) (1+\log (1+|\xi |^2 ))^2  R(t,\xi) (1+ R(t,\xi))}  
\end{align*}
for $|\xi | \geq \delta $, where $R(t,\xi) $ is given by \eqref{defR(t,xi)}. Thus, for $|\xi | \geq \sqrt{e-1} $ it holds that
$$  \left | \frac{b(\xi) -\sqrt{\log (1+|\xi |^2 ) }}{b(\xi)} \right |^2 \leq  \frac{ 1 }{ 15 \log ^2 (1+|\xi |^2 ) (1+\log (1+|\xi |^2 ))^4  }  .$$
Hence, one has
\begin{align*}
\int _{|\xi   | \geq \sqrt{e-1}} |  G_6(t, \xi) |^2 d \xi &= t^2 \int _{|\xi   | \geq \sqrt{e-1}} \left ( \frac{b(\xi) - \sqrt{\log(1+|\xi |^2)}}{b(\xi)} \right ) ^2 e^{-\frac{t}{1+ \log (1+|\xi |^2)}} \cos ^2 (\eta(\xi,t)) | \hat{u}_1 |^2 d \xi \\
&\leq \frac{1}{15}  t^2 \int _{|\xi   | \geq \sqrt{e-1}}  \frac{ e^{-\frac{t}{1+ \log (1+|\xi |^2)}} }{  \log ^2 (1+|\xi |^2 ) (1+\log (1+|\xi |^2 ))^4  }  \cos ^2 (\eta(\xi,t)) | \hat{u}_1 |^2 d \xi \\
&\leq \frac{4}{15}  t^2 \int _{|\xi   | \geq \sqrt{e-1}}  \frac{ e^{-\frac{t}{1+ \log (1+|\xi |^2)}} }{   (1+\log (1+|\xi |^2 ))^6  }  | \hat{u}_1 |^2 d \xi \\
&\leq C t^{-(l+4)} \|u_1 \|_{Y^l}^2 , \quad t>0.
\end{align*}

The estimates for $G_j(t,\xi) $ ($j=1,\cdots, 6$) together with the identity \eqref{u-perfil} provide the following result.
\begin{pro}\label{prophighfrequency} Let $n \geq 1 $, $l \geq 0 $ and $ (u_0, u_1) \in  Y^{l+1} \times Y^l  $. Then there exists a positive constant $C  $, which is independent of $t , u_0 $ and $u_1 $, such that 
\begin{align*} 
\int _{| \xi | \geq \sqrt{e-1}} & \left | \hat{u}(t, \xi)  - \varphi_2(t,\xi) \right |^2 d \xi \leq C \left ( \| u_0 \|_{Y^{l+1}}^2 + \| u_1 \| _{Y^l}^2 \right ) t^{-(l+3)}
\end{align*}
for $ t \geq 0 $, where $\varphi _2(t,\xi)$ is given by \eqref{defvarphi_2}. 
\end{pro}
\begin{flushright}
$\square$
\end{flushright}
\vspace{-0.7cm}


\subsection{Estimates on the whole space ${\bf R}^{n}$ }

In this subsection, we consider three special functions: $ \varphi _1(t,\xi) $, $\varphi _2(t,\xi)$, which are given by \eqref{defvarphi_1} and \eqref{defvarphi_2},  and  
\begin{align*}
 \varphi (\xi, t) &= \varphi _1 (t, \xi ) + \varphi _2(t,\xi)
\end{align*}
defined for $\xi \in {\bf R}_{\xi}^{n}$. 
\vspace{0.2cm}

We will prove that, under certain conditions, each of them is an asymptotic profile  as $t \rightarrow \infty $ of the solution $ \hat{u} (t,\xi) $ in ${\bf R}_{\xi}^{n}$.
\vspace{0.2cm}

\begin{lem}\label{lemmaperfilvarphi} Let $ n\geq 1$ and $ (u_0, u_1) \in (L^{1,1}({\bf R}^{n}) \cap Y^{l+1}) \times (L^{1,1}({\bf R}^{n}) \cap Y^{l})  $. Then there exists a constant $C>0 $, which is independent of $t, u_0, u_1$ such that 
\begin{align*}
\int _{{\bf R}^{n}} | \hat{u}(t,\xi) - \varphi(t,\xi) & |^2 d \xi \leq C\left ( t^{-\frac{n+2}{2}} + t^{-(l+3)} \right ) I_{0,l}^2 
\end{align*}
for $t \gg 1$, where
\begin{equation}\label{defiI_0}
I_{0,l} :=  \sqrt{\| u_0 \|_{1,1}^2 + \| u_1 \|_{1,1}^2 + \|u_0 \|_{Y^{l+1}}^2 + \| u_1 \|_{Y^l}^2 } .
\end{equation}
\end{lem} 
{\it Proof.}\
On the region $ | \xi | \leq \sqrt{e-1}$, the function $ \log (1+|\xi|^2 )   $ is positive and bounded by $1$, then for $t\geq0$ it holds that $$ \frac{-t}{\log (1+|\xi|^2 )} \leq -t . $$ 
We also have $ \sin a \leq a  $ for all $a \geq 0 $. Having this in mind we can get, for $t \geq 0 $, the estimates
\begin{align}\label{varphi_2baixafrequencia}
  \int _{|\xi | \leq \sqrt{e-1}} |\varphi _2 (t,\xi) |^2 d\xi & \leq 2 \int _{|\xi | \leq \sqrt{e-1} } e^{-\frac{t}{ \log (1+|\xi |^2 )}} \frac{\sin ^2 (\sqrt{\log (1+|\xi |^2 )}t )}{\log (1+|\xi |^2)} |\hat{u}_1|^2 d \xi \nonumber \\
  &+ 2 \int _{|\xi | \leq 1 } e^{-\frac{t}{ \log (1+|\xi |^2 )}} \cos ^2 (\sqrt{\log (1+|\xi |^2 )}t ) | \hat{u}_0| ^2 d \xi \nonumber \\
&\leq 2 t^2 e^{- t} \int _{|\xi | \leq \sqrt{e-1} }   |\hat{u}_1|^2 d \xi + 2 e^{- t}\int _{|\xi | \leq 1 }  | \hat{u}_0| ^2 d \xi \nonumber \\
&\leq 2 t^2 e^{- t} \| u_1 \|_2^2 + 2 e^{- t} \| u_0 \|_2^2 .
\end{align}


On the other hand, one knows that $$ e^{-2t \log (1+|\xi|^2) (1+ \log (1+|\xi|^2))}\leq e^{-2t \log(1+|\xi|^2) }  ,$$ because of  $ 1 + \log (1+|\xi |^2) \geq 1 $. Then
\begin{align*}
\int _{\eta \leq | \xi |}  |\varphi _1(t,\xi) |^2 d\xi &= |P_0 + P_1 |^2 \int _{| \xi | \geq \eta }e^{-2t\log(1+|\xi|^2)(1+\log(1+|\xi |^2))} d \xi \\
&\leq |P_0 + P_1 |^2 \int _{| \xi | \geq \eta } (1+|\xi |^2 )^{-2t} d \xi \\
&= |P_0 + P_1 |^2 \int _{ \eta \leq | \xi | \leq 1 } (1+|\xi |^2 )^{-2t} d \xi + |P_0 + P_1 |^2 \int _{| \xi | \geq 1 } (1+|\xi |^2 )^{-2t} d \xi \\
&=  \omega _n |P_0 + P_1 |^2 \int _{ \eta }^{1} (1+r^2 )^{-2t} r^{n-1}d r + \omega _n |P_0 + P_1 |^2 \int _{1 } ^{\infty}(1+r^2 )^{-2t} r^{n-1}d r \\
&\leq C |P_0 + P_1 |^2 \left (  (1+\eta ^2)^{-t} + \frac{2^{-t}}{t-1} \right )\\
& \leq  C \left ( \| u_0 \| _1^2 + \| u_1 \| _1^2 \right ) \left (  (1+\eta ^2)^{-t} + \frac{2^{-t}}{t-1} \right )
\end{align*}
with a generous constant $C > 0$ for $t \gg 1 $ due to Lemmas \ref{infit} and \ref{intermid}. We also note that both above estimates are of  exponential type.


\vspace{0.2cm}

Under these preparations we can get the desired estimate in the statement. At first, one notices that
\begin{align*}\label{4.23}
|\hat{u}(t,\xi) -\varphi (t,\xi)| &= |\hat{u}(t,\xi) -\varphi_1 (t,\xi) -\varphi_2 (t,\xi)| \leq |\hat{u}(t,\xi) -\varphi_1 (t,\xi)|  + |\varphi_2 (t,\xi)|.
\end{align*}
So, it holds that
\begin{align}
|\hat{u}(t,\xi) -\varphi (t,\xi)|^2  &\leq 2|\hat{u}(t,\xi) -\varphi_1 (t,\xi)|^2 + 2|\varphi_2 (t,\xi)|^2
\end{align}
Similarly, 
\begin{align}\label{4.24}
|\hat{u}(t,\xi) -\varphi (t,\xi)|^2 &\leq 2|\hat{u}(t,\xi) -\varphi_2 (t,\xi)|^2 + 2|\varphi_1 (t,\xi)|^2.
\end{align}
Also, one has 
\begin{align*}\label{4.25}
|\hat{u}(t,\xi) -\varphi (t,\xi)| \leq |\hat{u}(t,\xi)| +|\varphi_1 (t,\xi)|  + |\varphi_2 (t,\xi)|.
\end{align*}
By applying Young's inequality, we obtain 
\begin{align}
|\hat{u}(t,\xi) -\varphi (t,\xi)|^2 \leq 2|\hat{u}(t,\xi)|^2 +4|\varphi_1 (t,\xi)|^2  + 4|\varphi_2 (t,\xi)|^2.
\end{align}

Let us apply the estimates \eqref{4.23} on the region $|\xi |\leq \eta $, \eqref{4.24} on the region $|\xi|\geq \sqrt{e-1} $ and \eqref{4.25} on the middle frequency region $\eta\leq|\xi |\leq \sqrt{e-1}$, respectively. Then one can proceed the estimates as follows. 
\begin{align}
\int _{{\bf R}^{n}} |\hat{u}(t,\xi) - \varphi(t,\xi)|^2 d\xi &= \int _{|\xi |\leq \eta} |\hat{u}(t,\xi) - \varphi(t,\xi)|^2 d\xi + \int _{\eta\leq|\xi |\leq \sqrt{e-1} } |\hat{u}(t,\xi) - \varphi(t,\xi)|^2 d\xi \nonumber \\
&+  \int _{|\xi |\geq \sqrt{e-1} } |\hat{u}(t,\xi) - \varphi(t,\xi)|^2 d\xi \nonumber \\
&\leq 2 \int _{|\xi |\leq \eta}  |\hat{u}(t,\xi) - \varphi _1(t,\xi) |^2 d\xi + 2 \int _{|\xi |\leq \eta} |\varphi_2(t,\xi)|^2 d\xi \nonumber \\
&+ 2 \int _{\eta\leq|\xi |\leq \sqrt{e-1} } |\hat{u}(t,\xi)|^2 d\xi + 4 \int _{\eta\leq|\xi |\leq \sqrt{e-1}  } |\varphi_1(t,\xi)|^2 d\xi \nonumber \\
&+ 4 \int _{\eta\leq|\xi |\leq \sqrt{e-1} } |\varphi_2(t,\xi)|^2 d\xi  + 2 \int _{|\xi |\geq \sqrt{e-1} } |\hat{u}(t,\xi) - \varphi_2(t,\xi) |^2 d\xi \nonumber \\
&+ 2 \int _{|\xi |\geq \sqrt{e-1} } |\varphi_1(t,\xi)|^2 d\xi \nonumber \\
&\leq 2 \int _{|\xi |\leq \eta}  |\hat{u}(t,\xi) - \varphi _1(t,\xi) |^2 d\xi + 2 \int _{|\xi |\geq \sqrt{e-1} } |\hat{u}(t,\xi) - \varphi_2(t,\xi) |^2 d\xi \nonumber \\
&+ 4 \int _{|\xi |\leq \eta} |\varphi_2(t,\xi)|^2 d\xi + 4\int_{|\xi |\geq \sqrt{e-1} } |\varphi_1(t,\xi)|^2 d\xi \nonumber \\
&+ 2 \int _{\eta\leq|\xi |\leq \sqrt{e-1} } |\hat{u}(t,\xi)|^2 d\xi + 4 \int _{\eta\leq|\xi |\leq \sqrt{e-1}  } |\varphi_1(t,\xi)|^2 d\xi \nonumber \\
&+4\int _{\eta\leq|\xi |\leq \sqrt{e-1} } |\varphi_2(t,\xi)|^2 d\xi \nonumber \\
&=  2 \int _{|\xi |\leq \eta}  |\hat{u}(t,\xi) - \varphi _1(t,\xi) |^2 d\xi + 2 \int _{|\xi |\geq \sqrt{e-1} } |\hat{u}(t,\xi) - \varphi_2(t,\xi) |^2 d\xi \nonumber \\ 
&+ 4 \int _{|\xi |\leq \sqrt{e-1} } |\varphi_2(t,\xi)|^2 d\xi + 4\int_{|\xi |\geq \eta } |\varphi_1(t,\xi)|^2 d\xi \nonumber \\ 
&+2 \int _{\eta\leq|\xi |\leq \sqrt{e-1} } |\hat{u}(t,\xi)|^2 d\xi. 
\end{align}

Furthermore, Propositions 4.1 and 4.2 tell us that 
\begin{align*}
&\int _{|\xi |\leq \eta } |\hat{u}(t,\xi) - \varphi _1(t,\xi)|^2 d\xi \leq C I_0^2 t^{-\frac{n+2}{2}}, \\
&\int _{|\xi| \geq \sqrt{e-1}} |\hat{u}(t,\xi)- \varphi _2(t,\xi)|^2 d\xi \leq CI_0^2 t^{-(l+3)} 
\end{align*}
for $t \gg1 $. By combining the previous two estimates with Propositions \ref{proplowfrequency} and \ref{prophighfrequency} and Lemma \ref{lemmadecayulow-middle} we can conclude the desired estimate 
$$\int _{{\bf R}^{n}} |\hat{u}(t,\xi) - \varphi(t,\xi)|^2 d\xi \leq \tilde{C}I_0^2 (t^{-\frac{n+2}{2}}+t^{-(l+3)}) \quad t \gg 1. $$

\hfill
$\Box$
\vspace{0.2cm}

\begin{lem}\label{decayvarphi_1} Let $u_0, u_1 \in L^{1}({\bf R}^{n}) $ and the function $\varphi _1 (t,\xi) $ is defined in \eqref{defvarphi_1}. Then there exists positive constants $ C_1, C_2$, depending only on dimension $n $, such that 
\begin{equation}
C_1 |P_0+P_1|^2 t^{-\frac{n}{2}}  \leq \int_{{\bf R}_{\xi}^{n}} |\varphi _1 (t,\xi)|^2 d\xi \leq C_2 (\|u_0\|_1^2 +\|u_1\|_1^2) t^{-\frac{n}{2}}
\end{equation}
for $t \gg 1 $.
\end{lem}
{\it Proof.}\
The function $ \varphi _1(t, \xi)  $ satisfies 
$$ | \varphi _1(t,\xi ) | \leq |P_0 + P_1| (1+|\xi |^2) ^{-t} $$
for $ \xi \in {\bf R}_{\xi}^{n}$, since $ 1+ \log (1+|\xi|^2) \geq 1. $
By using Lemmas \ref{general-p} and \ref{infit}, we immediately concluded that
\begin{align*}
 \int _{{\bf R}_{\xi}^{n}} |\varphi _1 (t, \xi)|^2 d\xi &\leq |P_0+P_1|^2 \int _{{\bf R}_{\xi}^{n}} (1+|\xi |^2)^{-2t} d \xi \\
 &=  \omega _n |P_0+P_1|^2 \int _{0}^{1} (1+r^2)^{-2t} r^{n-1}d r + \omega _n |P_0+P_1|^2 \int _{1}^{\infty} (1+r^2)^{-2t} r^{n-1}d r \\
 &\leq C \omega _n  |P_0+P_1|^2 \left ( t^{-\frac{n}{2}} + \frac{2^{-t}}{t-1} \right ) \\
 &\leq  C \omega _n  (\|u_0\|_1^2 +\|u_1\|_1^2) t^{-\frac{n}{2}}
\end{align*}
for $ t \gg 1$.

\vspace{0.2cm}
On the other hand, for $| \xi | \leq \eta $, we have $ 1+ \log (1+|\xi |^2)  \leq 1+\log (1+|\eta|^2)=k_\eta $. Thus, $$  | \varphi(t,\xi ) | \geq  |P_0 + P_1| (1+|\xi |^2) ^{- k_\eta t}  $$
for $| \xi | \leq \eta $. First, we choose $ t_0>0 $ such that, for all $t> t_0  $ it holds that $ t^{- \frac{1}{2}} \leq \eta  $, and 
$$   \frac{1}{e^{4 k_\eta }}\leq \left ( 1+ \frac{1}{t} \right )^{-2k_\eta t} \leq 1.$$  
Such $t_0$ exists, because one has
$$ \lim _{t \rightarrow \infty } \left ( 1+ \frac{1}{t} \right )^{-2k_\eta t} = \frac{1}{e^{2k_\eta}}.$$
For this choice, we can compute as follows:
\begin{align*}
\int _{{\bf R}_{\xi}^{n}} |\varphi _1(t,\xi) |^2 d \xi & \geq \int _{|\xi | \leq \eta } |\varphi _1(t,\xi) |^2 d \xi \geq \omega _n |P_0+P_1|^2 \int _{0}^{\eta } (1+r^2)^{-2k_\eta t} r^{n-1}dr \\
&\geq \omega _n |P_0+P_1|^2 \int _{0}^{t^{-\frac{1}{2}} } (1+r^2)^{-2k_\eta t} r^{n-1}dr   \\ 
&\geq  \omega _n |P_0+P_1|^2 \left ( 1 +\frac{1}{t} \right )^{-2k_{\eta} t} \int _{0}^{t^{-\frac{1}{2}} }  r^{n-1}dr   \\ 
&=  \frac{\omega _n}{n} |P_0+P_1|^2 \left ( 1 +\frac{1}{t} \right )^{-2k_{\eta} t} t^{-\frac{n}{2}}\\
&\geq \frac{\omega _n e^{-4k_{\eta}}}{n}|P_0+P_1|^2 t^{-\frac{n}{2}}
\end{align*}
for $ t > t_0$.

\hfill
$\Box$
\vspace{0.2cm}

\begin{lem}\label{decayvarphi_2} Let $n\geq 1 $, $ l \geq 0$ and $ (u_0, u_1) \in Y^{l+1} \times Y^l $. Then there exists a constant $ C>0$, which is independent   of $ u_0, u_1 $ and $t$, such that
$$ \int _{{\bf R}_{\xi}^{n}} |\varphi _2(t, \xi)|^2 d\xi \leq C I_{0,l}^2 t^{-(l+1)}  $$
for all $t>0$, where $I_{0,l}$ is given in \eqref{defiI_0}.
\end{lem}
{\it Proof.}\
By definition of $\varphi _2(t,\xi) $ in \eqref{defvarphi_2}, we have
\begin{align*}
|\varphi _2(t, \xi) | \leq e^{-\frac{t}{2 \log (1+|\xi |^2 )}} \frac{|\sin (\sqrt{\log (1+|\xi |^2 )}t )|}{\sqrt{\log (1+|\xi |^2)}} |\hat{u}_1| +  e^{-\frac{t}{2 \log (1+|\xi |^2 )}} |\cos (\sqrt{\log (1+|\xi |^2 )}t )| |\hat{u}_0|.
\end{align*}
Hence, the Young's inequality enable us to get 
\begin{align}\label{eqvarphi_2^2}
|\varphi _2(t, \xi) |^2 \leq 2e^{-\frac{t}{\log (1+|\xi |^2 )}} \frac{\sin ^2(\sqrt{\log (1+|\xi |^2 )}t )}{\log (1+|\xi |^2)} |\hat{u}_1|^2 +  2e^{-\frac{t}{ \log (1+|\xi |^2 )}} \cos^2 (\sqrt{\log (1+|\xi |^2 )}t )|\hat{u}_0|^2.
\end{align}
It follows from \eqref{varphi_2baixafrequencia}, we get
\begin{align}\label{estimatevarphi_2low}
\int _{|\xi| \leq \sqrt{e-1}} |\varphi _2(t, \xi)|^2 d\xi &\leq 2 t^2 e^{-t} \|u_1\|_2^2 + 2e^{-t} \|u_0\|_2^2, \quad t>0. 
\end{align}

On the high frequency zone $ | \xi | \geq \sqrt{e-1} $ it holds that
\begin{align}\label{eq4.24}
\frac{1}{1+\log(1+|\xi|^2) } \leq \frac{1}{\log(1+|\xi|^2) } \leq  \frac{2}{1+ \log(1+|\xi|^2 )}.
\end{align}
By using the inequality \eqref{eq4.24} and the estimates \eqref{decay-nu} and \eqref{eqvarphi_2^2}, one can obtain
\begin{align}\label{estimatevarphi_2high}
\int _{|\xi| \geq \sqrt{e-1}} |\varphi _2(t, \xi)|^2 d\xi &\leq 2 \int _{|\xi| \geq \sqrt{e-1}} e^{-\frac{t}{ \log (1+|\xi |^2 )}} \frac{\sin^2 (\sqrt{\log (1+|\xi |^2 )}t )}{\log (1+|\xi |^2) } |\hat{u}_1|^2 d \xi \nonumber \\
&+ 2 \int _{|\xi| \geq \sqrt{e-1}}  e^{-\frac{t}{ \log (1+|\xi |^2 )}} \cos^2  (\sqrt{\log (1+|\xi |^2 )}t ) |\hat{u}_0|^2 d \xi \nonumber \\
&\leq 4 \int _{|\xi| \geq \sqrt{e-1}} e^{-\frac{t}{1+ \log (1+|\xi |^2 )}} \frac{1}{1+\log (1+|\xi |^2) } |\hat{u}_1|^2 d \xi  \nonumber \\
&+ 2 \int _{|\xi| \geq \sqrt{e-1}}  e^{-\frac{t}{ 1+ \log (1+|\xi |^2 )}}  |\hat{u}_0|^2 d \xi \nonumber \\
&=  4 \int _{|\xi| \geq \sqrt{e-1}}  \frac{e^{-\frac{t}{1+ \log (1+|\xi |^2 )}}}{(1+\log (1+|\xi |^2))^{l+1} } (1+\log (1+|\xi |^2))^{l} |\hat{u}_1|^2 d \xi \nonumber \\
&+ 2 \int _{|\xi| \geq \sqrt{e-1}}  \frac{e^{-\frac{t}{1+ \log (1+|\xi |^2 )}}}{(1+\log (1+|\xi |^2))^{l+1} } (1+\log (1+|\xi |^2))^{l+1} |\hat{u}_0|^2 d \xi \nonumber \\
&\leq 4 C t^{-(l+1)} \|u_1 \|_{Y^{l}}^2 + 2 C t^{-(l+1)}\|u_0 \|_{Y^{l+1}}^2, \quad t>0.
\end{align}
By combining estimates \eqref{estimatevarphi_2low} and \eqref{estimatevarphi_2high} one can conclude 
$$ \int _{{\bf R}_{\xi}^{n}} |\varphi _2(t, \xi)|^2 d\xi \leq C t^{-(l+1)} ( \|u_1 \|_{Y^{l}}^2 +\|u_0 \|_{Y^{l+1}}^2 ),\quad t \gg 1$$
for some generous constant $C>0 $, independent of $ u_0, u_1 $ and $t$. 
\hfill
$\Box$

\subsection{Proof of Theorems \ref{1.1}, \ref{1.2} and \ref{1.3}}

Now, let us prove our main Theorems \ref{1.1}, \ref{1.2}, and \ref{1.3} at a stroke in the following paragraph.\\

First, it follows from Lemma \ref{lemmaperfilvarphi} we have the estimate 
 \begin{equation}\label{28}
 \int _{{\bf R}_{\xi}^{n}} | \hat{u}(t,\xi) - \varphi(t,\xi)  |^2 d \xi  \leq C I_0^2 \left ( t^{-\frac{n+2}{2}} + t^{-(l+3)} \right ) =:P_n(t) , 
 \end{equation}
 where $ I_0 $ is defined by \eqref{defiI_0}. 
 \vspace{0.2cm}
 
Next, since we can write as $ \hat{u}(t,\xi) - \varphi _1 (t,\xi) = \hat{u}(t,\xi) - \varphi (t,\xi) + \varphi_2(t,\xi) $, and $ \hat{u}(t,\xi) - \varphi _2 (t,\xi) = \hat{u}(t,\xi) - \varphi (t,\xi) + \varphi_1(t,\xi)  $, from Lemmas \ref{decayvarphi_2} and \ref{decayvarphi_1} one has
\[\int _{{\bf R}_{\xi}^{n}} | \hat{u}(t,\xi) - \varphi _1 (t,\xi) | ^2 d\xi \leq 2 \int _{{\bf R}_{\xi}^{n}} | \hat{u}(t,\xi) - \varphi (t,\xi) | ^2 d\xi + 2 \int _{{\bf R}_{\xi}^{n}} |  \varphi _2 (t,\xi) | ^2 d\xi \]
\begin{equation}\label{29}
\leq 2C I_0^2 \left ( t^{-\frac{n+2}{2}} + t^{-(l+3)} + t^{-(l+1)} \right ) \leq 4C I_0^2 \left ( t^{-\frac{n+2}{2}}  + t^{-(l+1)} \right ) =:M_n(t),
\end{equation}
 and 

\[\int _{{\bf R}_{\xi}^{n}} | \hat{u}(t,\xi) - \varphi _2 (t,\xi) | ^2 d\xi \leq 2 \int _{{\bf R}_{\xi}^{n}} | \hat{u}(t,\xi) - \varphi (t,\xi) | ^2 d\xi + 2 \int _{{\bf R}_{\xi}^{n}} |  \varphi _1 (t,\xi) | ^2 d\xi\]
 \begin{equation}\label{30}
\leq 2C I_0^2 \left ( t^{-\frac{n+2}{2}} + t^{-(l+3)} + t^{-\frac{n}{2} } \right ) \leq 4C I_0^2 \left (  t^{-(l+3)} + t^{-\frac{n}{2} } \right ) =: Q_n(t).
 \end{equation}
By an asymptotic profile we mean the part of $\hat{u}(t,\xi)$ that decays with the slowest time rate. According to Lemmas \ref{decayvarphi_1} and \ref{decayvarphi_2} we know that 
$$\int _{{\bf R}_{\xi}^{n}} | \varphi _1(t,\xi) |^2 \leq CI_0^2 t^{-\frac{n}{2}}, \quad \int _{{\bf R}_{\xi}^{n}} | \varphi _2(t,\xi) |^2 \leq CI_0^2 t^{-(l+1)}, \quad \int _{{\bf R}_{\xi}^{n}} | \varphi (t,\xi) |^2 \leq CI_0^2 (t^{-\frac{n}{2}} +t^{-(l+1)} ).$$
Therefore, the asymptotic profile of the solution $\hat{u}(t,\xi)$ as $t \to +\infty$ is 
\begin{itemize}
\item[{\rm (i)}. ] $ \varphi _1(t,\xi) $ if $ \frac{n}{2}< l+1 $ (compare with \eqref{29}), 
\item[{\rm (ii)}. ] $ \varphi (t,\xi) $ if $ \frac{n}{2} = l+1 $ (compare with \eqref{28}),
\item[{\rm (iii)}. ] $ \varphi _2(t,\xi) $ if $ \frac{n}{2} > l+1 $ (compare with \eqref{30}).
\end{itemize}


In the final part of this subsection we will discuss the decay rate of $ P_n(t), M_n(t) $ and $ Q_n(t) $ related to the differences between the solution $\hat{u}(t,\xi) $ and the suitable asymptotic profiles.\\
Let $l \geq 1$ as this is necessary for the existence and uniqueness of the solution (see Theorem 2.2). For this purpose we introduce a (critical) value $l^{*}(n)$ on the regularity $l \geq 1$ of the initial data such that
\[l^{*}(n) := \frac{n}{2}-1.\]
Let us prove Theorems \ref{1.1}, \ref{1.2} and \ref{1.3} as follows.
 
\begin{itemize}
\item[{\rm (i)}. ] The case for $l^{*}(n) < l$. 
\\
First, we consider $l=1$. In this case, $\varphi _1(t,\xi) $ is asymptotic profile for $ n < 4 $. 

$ \bullet $ If $ n \leq 2 $, then  $ \frac{n+2}{2} \leq 2=l+1$, and so $M_n(t)\leq 8 C I_0^2 t^{-\frac{n+2}{2}}$. \\
$ \bullet $ If $ n = 3  $, then $ 2= l+1 < \frac{n+2}{2}  =  \frac{5}{2}   $. Thus $M_n(t)\leq 8 C I_0^2 t^{-(l+1)} = 8 C I_0^2 t^{-2}$. 

Now let us consider the case $ l>1 $. \\
$ \bullet $ In this case, the rate $ t^{-(l+1)} $ is better than $ t^{-2} $. Therefore, if $ n \leq 2 $, we have $ M_n(t)\leq 8 C I_0^2 t^{-\frac{n+2}{2}}$. \\
$ \bullet $ If  $ n >2  $ and $ \frac{n}{2} \leq l  $, then $ l > 1 $ and  $ M_n(t)\leq 8 C I_0^2 t^{-\frac{n+2}{2}}$. \\
$\bullet $ If $ n \geq 4  $ and $ \frac{n}{2}-1 < l < \frac{n}{2} $, we obtain $ l>1  $ and $ M_n (t) \leq 8 C I_0^2 t^{-(l+1)}$. \\
$ \bullet $ For $ n =3 $, we need $1< l \leq \frac{3}{2} $, in order that $ M_n (t) \leq 8 C I_0^2 t^{-(l+1)}$.\\
These observations together with the Plancherel Theorem imply the desired statement of Theorem \ref{1.1}.

\item[{\rm (ii)}. ] The case for $l^{*}(n) > l$.\\
If $ l =1  $, $\varphi _2 (t,\xi) $ is asymptotic profile for $ n>4 $.\\
$\bullet $ If $ 4< n\leq 8  $, then $\frac{n}{2} \leq 4 =l+3 $. So $ Q_n(t) \leq 8C I_0^2 t^{- \frac{n}{2}} $. \\
$ \bullet $ For  $n>8 $, we have $ \frac{n}{2}>4=l+3 $ and $  Q_n(t) \leq  8C I_0^2 t^{-(l+3)} = 8C I_0^2 t^{-4} $.\\
By assuming $ l>1 $, it is necessary that $n>2l+2>4 $.\\
$ \bullet $ If $ 4<n\leq 8 $, then $ \frac{n}{2} \leq 4 < l+3$. Therefore,  $  Q_n(t) \leq 8C I_0^2 t^{- \frac{n}{2}} $. \\
$\bullet $ For $ n> 8 $ and $ \frac{n}{2}-3 < l < \frac{n}{2} -1$, we have $ l>1 $ and $  Q_n(t) \leq 8C I_0^2 t^{- \frac{n}{2}} $.\\
$\bullet $  If $n>8 $ and $ 1<l\leq \frac{n}{2}-3 $, we obtain $  Q_n(t) \leq 8C I_0^2 t^{-(l+3)} $.\\
These observations together with the Plancherel Theorem imply the desired statement of Theorem \ref{1.2}.

\item[{\rm (iii)}. ] The case for $l^{*}(n) = l$. \\
This condition implies that $ \frac{n+2}{2} = l+2 < l+3 $. Then we have $P_n(t) \leq 2 C I_0^2 t^{-\frac{n+2}{2}} $. Due to $ l\geq 1 $, this estimate holds only for $ n \geq 4 $. \\
These observations together with the Plancherel Theorem imply the desired statement of Theorem \ref{1.3}.

\end{itemize} 




\subsection{Optimal decay rates of the solution}

In this section we prove Theorems \ref{TEO0}, \ref{TEO1} and \ref{TEO2} based on previously obtained decay estimates.  
\vspace{0.2cm}  

We have already used the decomposition such as
$$ \hat{u}(t,\xi) = \hat{u}(t,\xi) - \varphi (t,\xi) + \varphi _1(t,\xi) +\varphi_2(t,\xi), $$
where $ \varphi (t,\xi) = \varphi _1(t,\xi) + \varphi_2(t,\xi) $ with $ \varphi _1(t,\xi) $ and $ \varphi _2 (t,\xi)$ are given by \eqref{defvarphi_1}, \eqref{defvarphi_2}. Since $ u_0 $ and $u_1 $ have the required regularity in Lemmas \ref{lemmaperfilvarphi}, \ref{decayvarphi_1} and \ref{decayvarphi_2}, one can get 
 \begin{align}\label{expressaou}
 \int _{{\bf R}_{\xi}^{n}} |\hat{u}(t,\xi)|^2 dx &\leq 4  \int _{{\bf R}_{\xi}^{n}} | \hat{u}(t,\xi) -  \varphi (t,\xi)|^2 d\xi + 4 \int _{{\bf R}_{\xi}^{n}} |\varphi _1(t,\xi)|^2 d\xi +4 \int _{{\bf R}_{\xi}^{n}} |\varphi_2(t,\xi)|^2 d\xi \nonumber \\
&\leq 4C I_0^2 (t^{-\frac{n+2}{2}}+t^{-(l+3)} + t^{-\frac{n}{2}} + t^{-(l+1)}) \quad (t \gg 1)\nonumber \\
&\leq K I_0^2  ( t^{-\frac{n}{2}} + t^{-(l+1)})=:R_n(t) 
 \end{align}
with some constant $K>0$. 
 
In the same way as we did in the previous subsection, we compare $\displaystyle{\frac{n}{2}}$ and $(l+1)$ in order to obtain the decay rate of the solution. 
\begin{itemize}
\item[{\rm (i)}. ] If $ l \geq 1 $ and $ n \leq 3  $, then $\displaystyle{\frac{n}{2}} < 2 \leq l+1 $. So $R_n(t) \leq 2 K I_0^2 t^{-\frac{n}{2}}$.
\item[{\rm (ii)}. ] If $ n>4 $ and $ l> \displaystyle{\frac{n}{2}} -1 $, we have  $R_n(t) \leq 2 K I_0^2 t^{-\frac{n}{2}}$.
\item[{\rm (iii)}. ]If $ n\geq 4 $ and $ 1 \leq l \leq \displaystyle{\frac{n}{2}}-1 $, then $ l+1 \leq \displaystyle{\frac{n}{2}}$. Thus $ R_n(t) \leq 2 K I_0^2 t^{-(l+1)} $.
\end{itemize}
The last item (iii) combined with the expression \eqref{expressaou} and Plancherel Theorem completes the proof of Theorem \ref{TEO0}.
\hfill
$\Box$

\vspace{0.2cm}
\par
Items (i) and (ii) give us conditions for the decay rate of the solution to be better than $t^{-\frac{n}{2}} $. Furthermore, we may prove that this rate is optimal under these same conditions. Recalling the fact that the condition $ (u_0, u_1) \in Y^2 \times Y^1  $ is necessary for the existence and uniqueness of the solution $u(t,x)$ to problem (1.1)-(1.2) according to Theorem \ref{existence}, we can prove Theorems \ref{TEO1} and \ref{TEO2} at a stroke.\\



{\it Proof of Theorems \ref{TEO1} and \ref{TEO2}.}\
\vspace{0.2cm}

We first observe the case of $ (u_0, u_1) \in Y^{l+1} \times Y^l $ with $ l > \displaystyle{\frac{n-2}{2}} $.
On the one hand, 
\begin{align*}
| \hat{u}(t,\xi) |^2 &\leq \left ( |\hat{u}(t,\xi) -\varphi (t,\xi) | + |\varphi(t,\xi)| \right )^2 \\
&\leq 2|\hat{u}(t,\xi) -\varphi (t,\xi) |^2 + 2|\varphi(t,\xi)|^2 \\
&\leq 4|\hat{u}(t,\xi) -\varphi (t,\xi) |^2 + 4|\varphi_1(t,\xi)|^2 + 4|\varphi_2(t,\xi)|^2 .
\end{align*}
So, 
\begin{align*}
\int _{{\bf R}_{\xi}^{n}} | \hat{u}(t,\xi) |^2 d\xi &\leq 4 \int _{{\bf R}_{\xi}^{n}} |\hat{u}(t,\xi) -\varphi (t,\xi) |^2 d\xi + 4\int _{{\bf R}_{\xi}^{n}} |\varphi_1(t,\xi)|^2 d\xi + 4\int _{{\bf R}_{\xi}^{n}} |\varphi_2(t,\xi)|^2 d\xi \\
&\leq   CI_{0,l}^2 \left ( t^{- \frac{n+2}{2}} + t^{-\frac{n}{2}} + t^{-(l+1)} \right ) \\
&\leq CI_{0,l}^2 t^{-\frac{n}{2}},
\end{align*}
since $ l+ 1 > \displaystyle{\frac{n}{2}} $, due to Lemmas \ref{decayvarphi_1}, \ref{decayvarphi_2} and \ref{lemmaperfilvarphi}. Thus, the upper bound estimates in Theorems \ref{TEO1} and \ref{TEO2} can be proved by choosing $l = 1$ and $l = \displaystyle{\frac{n-2}{2}} + \varepsilon$, respectively.\\

On the other hand, since one has $ |\varphi (t,\xi) |\leq  |\varphi (t,\xi) - \hat{u}(t,\xi)| + |\hat{u}(t,\xi)|  $ and $| \varphi _1 (t,\xi)| \leq |\varphi _1 (t,\xi) + \varphi _2 (t,\xi) |+ |\varphi _2(t,\xi)|$, by using  Young's inequality, we obtain 
\begin{align*}
 |\hat{u}(t,\xi)|^2 &\geq \frac{1}{2} |\varphi  (t,\xi) |^2 - |\varphi (t,\xi) - \hat{u}(t,\xi)|^2 \\
&\geq \frac{1}{4} | \varphi _1 (t,\xi)|^2 -  \frac{1}{2}|\varphi _2(t,\xi)|^2 - |\varphi (t,\xi) - \hat{u}(t,\xi)|^2 .
\end{align*} 
By using the estimates just obtained in Lemmas \ref{decayvarphi_1}, \ref{decayvarphi_2} and \ref{lemmaperfilvarphi}, one can obtain the following expression:
\begin{align}\label{eq5.3}
\int _{{\bf R}_{\xi}^{n}} |\hat{u}(t,\xi)|^2 d \xi &\geq \frac{1}{4} \int _{{\bf R}_{\xi}^{n}} |\varphi _1 (t,\xi) |^2 d \xi - \frac{1}{2} \int _{{\bf R}_{\xi}^{n}} |\varphi _2 (t,\xi) |^2 d \xi - \int_{{\bf R}_{\xi}^{n}} | \hat{u}(t,\xi) - \varphi (t,\xi) |^2 d\xi \nonumber \\
&\geq C_1 |P_0+P_1|^2 t^{-\frac{n}{2}} - C I_{0,l}^2 t^{-(l+1)}  - C I_{0,l}^2 t^{-\frac{n+2}{2}} - C I_{0,l}^2 t^{-(l+3)} \nonumber \\
&= t^{-\frac{n}{2}} \left ( C_1 |P_0+P_1|^2 - CI_0^2 t^{-\frac{2l-n+2}{2} } - CI_0^2 t^{-1} - CI_0^2 t^{-\frac{2l-n+6}{2} } \right ).
\end{align}
If $ \frac{n}{2} < l+1  $,  then $ \frac{2l-n+2}{2} >0 $ and $ \frac{2l-n+6}{2} >0 $, because of $ l+1 < l+3 $. Hence, one has
$$ \lim _{t \rightarrow \infty} \left (  CI_0^2 t^{-\frac{2l-n+6}{2} }+ CI_0^2 t^{-1} + CI_0^2 t^{-\frac{2l-n+2}{2} } \right ) = 0,$$
so that there exists $ t_1 \gg 1$ such that 
$$  CI_0^2 t^{-\frac{2l-n+6}{2} }+ CI_0^2 t^{-1} + CI_0^2 t^{-\frac{2l-n+2}{2} } \leq \frac{C_1}{2}|P_0+P_1|^2
$$
 for all $ t \geq t_1 $ in the case of $\vert P_{1} + P_{0}\vert \ne 0$. That is, for $ t\geq t_1$ it holds that 
 $$ C_1|P_0+P_1|^2 - CI_0^2 t^{-\frac{2l-n+6}{2} } - CI_0^2 t^{-1} - CI_0^2 t^{-\frac{2l-n+2}{2} } \geq \frac{C_1}{2}|P_0+P_1|^2. $$
 Therefore, one can arrive at the crucial estimate
 \begin{align}
 \int _{{\bf R}_{\xi}^{n}} |\hat{u}(t,\xi)|^2 d \xi &\geq  \frac{C_1}{2}|P_0+P_1|^2  t^{-\frac{n}{2}} 
 \end{align}
for $ t \geq t_1 $ because of \eqref{eq5.3}. By choosing $l = 1$ in Theorem \ref{TEO1}, and $l = \displaystyle{\frac{n-2}{2}} + \varepsilon$ in Theorem \ref{TEO2}, one can get the desired estimates.
\hfill
$\Box$
\vspace{0.2cm}

\par
\vspace{0.5cm}
\noindent{\em Acknowledgement.}
\smallskip
The work of the first author (R. C. CHAR\~AO) was partially supported by PRINT/ CAPES - Process 88881.310536/2018-00 and the work of the third author (R. IKEHATA) was supported in part by Grant-in-Aid for Scientific Research (C)20K03682  of JSPS. 




\end{document}